%% file: main.tex
\documentclass[onefignum,onetabnum]{siamonline190516}

\input{packages}
\input{macros}

\headers{Parallelism Approaches for GCP with GenTen}{J. M. Myers And
E. T. Phipps}

\title{Synchronous and Asynchronous Parallelism Approaches for
  Generalized Canonical
  Polyadic Tensor Decomposition with GenTen\thanks{Sandia National
    Laboratories is
    a multimission laboratory managed and operated by National
    Technology \& Engineering Solutions of Sandia, LLC, a wholly owned
    subsidiary of Honeywell International Inc., for the U.S. Department of
    Energy’s National Nuclear Security Administration under contract
DE-NA0003525. SAND2026-20026O.}}

\author{Jeremy M. Myers\thanks{Sandia National Laboratories
    (\email{jermyer@sandia.gov},
\email{etphipp@sandia.gov})} \and Eric T. Phipps\footnotemark[2] }

\ifpdf
\hypersetup{ pdftitle={Synchronous and Asynchronous Parallelism
    Approaches for Generalized Canonical
  Polyadic Tensor Decomposition with GenTen}, pdfauthor={J.
M. Myers and E. T. Phipps} }
\fi

\begin{document}

\maketitle

\begin{abstract}
  \input{abstract}
\end{abstract}




\section{Introduction}
\label{intro}
\input{intro}

\section{Background}
\label{background}
\input{background}

\section{Hybrid Parallel GCP with GenTen}
\label{genten-dist-par}
\input{genten-distributed-parallelism}

\section{GCP-FedAdam}
\label{fed-opt}
\input{fed-opt}

\section{Numerical experiments}
\label{experiments}
\input{experiments}

\section{Conclusions}
\label{conc}
\input{conclusions}

\clearpage
\section*{Acknowledgments}
We thank Theresa Portone of Sandia National Laboratories for her
expertise is Latin hypercube sampling, using Dakota, and interpreting
simulation data.
This material is based upon work partially supported by the U.S. Department of Energy, Office of Science,
Office of Advanced Scientific Computing Research under contract number DE-NA0003525.

\clearpage
\bibliographystyle{siamplain}
\bibliography{refs}

\end{document}

%% file: packages.tex
\usepackage{csquotes}
\usepackage{bibentry}
\usepackage{cite}
\usepackage{epigraph}

\usepackage{bm}
\usepackage[mathscr]{euscript}
\usepackage{stmaryrd}
\SetSymbolFont{stmry}{bold}{U}{stmry}{m}{n} 

\usepackage{hyperref}
\usepackage{cleveref}

\usepackage{caption}
\usepackage{float}
\usepackage{graphicx}
\graphicspath{{./images/}}
\usepackage{subcaption}
\usepackage{tikz}
\usetikzlibrary{calc}
\usepackage{wrapfig}
\usepackage{pgfplots}
\pgfplotsset{compat=1.18}
\pgfplotsset{select coords between index/.style 2 args={ x filter/.code={
      \ifnum\coordindex<#1\fi
\ifnum\coordindex>#2\fi } }}
\usepackage{booktabs}
\usepackage{tabularx}

\usepackage{adjustbox}

\usepackage{amsfonts,amsmath,amssymb,amsbsy}
\usepackage{mathtools}
\usepackage{etoolbox}

\usepackage{algorithm}
\usepackage[noend]{algpseudocode}
\Crefname{ALC@unique}{Line}{Lines} 
\algnewcommand{\LeftComment}[1]{\{#1\}\hfill}

\usepackage{enumitem}
\usepackage{siunitx}
\usepackage{url}

\usepackage{xparse}

\newcommand{\TODO}[2]{}
\newcommand{\EnableTODO}{%
\renewcommand{\TODO}[2]{\textcolor{red}{\textbf{\textit{[##1---##2]}}}}}
 \EnableTODO

%% file: macros.tex
\NewDocumentCommand{\N}{o}{\mathbb{N}}
\NewDocumentCommand{\R}{o}{\mathbb{R}}
\NewDocumentCommand{\F}{o}{\mathbb{F}}
\NewDocumentCommand{\Z}{o}{\mathbb{Z}}

\NewDocumentCommand{\BigO}{o}{\mathcal{O}}

\makeatletter
\newcommand{\@giventhatstar}[2]{\left(#1\;\middle|\;#2\right)}
\newcommand{\@giventhatnostar}[3][]{#1(#2\;#1|\;#3#1)}
\newcommand{\giventhat}{\@ifstar\@giventhatstar\@giventhatnostar}
\makeatother

\NewDocumentCommand{\tns}{m}{\boldsymbol{\mathscr{#1}}}

\NewDocumentCommand{\mat}{m}{\boldsymbol{\mathbf{#1}}}

\NewDocumentCommand{\mind}{o}{\boldsymbol{\mathbf{i}}}

\DeclarePairedDelimiter\abs{\lvert}{\rvert}%
\DeclarePairedDelimiter\norm{\lVert}{\rVert}%

\NewDocumentCommand{\gcp}{}{\textsc{GCP}}
\NewDocumentCommand{\adam}{}{\textsc{Adam}}
\NewDocumentCommand{\sgd}{}{\textsc{SGD}}

\NewDocumentCommand{\fedopt}{}{\textsc{FedOpt}}
\NewDocumentCommand{\fedadam}{}{\textsc{FedAdam}}
\NewDocumentCommand{\gcpsgd}{}{{\gcp}-{\sgd}}

\NewDocumentCommand{\gcpfed}{}{{\gcp}-{\fedadam}}

\makeatletter
\let\oldabs\abs \def\abs{\@ifstar{\oldabs}{\oldabs*}}
\let\oldnorm\norm \def\norm{\@ifstar{\oldnorm}{\oldnorm*}}
\makeatother

\makeatletter
\let\save@mathaccent\mathaccent \newcommand*\if@single[3]{%
  \setbox0\hbox{${\mathaccent"0362{#1}}^H$}%
  \setbox2\hbox{${\mathaccent"0362{\kern0pt#1}}^H$}%
\ifdim\ht0=\ht2 #3\else #2\fi }
\newcommand*\rel@kern[1]{\kern#1\dimexpr\macc@kerna}
\newcommand*\widebar[1]{\@ifnextchar^{{\wide@bar{#1}{0}}}{\wide@bar{#1}{1}}}
\newcommand*\wide@bar[2]{\if@single{#1}{\wide@bar@{#1}{#2}{1}}{\wide@bar@{#1}{#2}{2}}}
\newcommand*\wide@bar@[3]{%
  \begingroup
  \def\mathaccent##1##2{%
    \let\mathaccent\save@mathaccent
    \if#32 \let\macc@nucleus\first@char \fi
    \setbox\z@\hbox{$\macc@style{\macc@nucleus}_{}$}%
    \setbox\tw@\hbox{$\macc@style{\macc@nucleus}{}_{}$}%
    \dimen@\wd\tw@ \advance\dimen@-\wd\z@
    \divide\dimen@ 3 \@tempdima\wd\tw@ \advance\@tempdima-\scriptspace
    \divide\@tempdima 10 \advance\dimen@-\@tempdima
    \ifdim\dimen@>\z@ \dimen@0pt\fi
    \rel@kern{0.6}\kern-\dimen@ \if#31
    \overline{\rel@kern{-0.6}\kern\dimen@\macc@nucleus\rel@kern{0.4}\kern\dimen@}%
    \advance\dimen@0.4\dimexpr\macc@kerna
    \let\final@kern#2%
    \ifdim\dimen@<\z@ \let\final@kern1\fi \if\final@kern1 \kern-\dimen@\fi
    \else
    \overline{\rel@kern{-0.6}\kern\dimen@#1}%
    \fi
  }%
  \macc@depth\@ne \let\math@bgroup\@empty \let\math@egroup\macc@set@skewchar
  \mathsurround\z@ \frozen@everymath{\mathgroup\macc@group\relax}%
  \macc@set@skewchar\relax \let\mathaccentV\macc@nested@a
  \if#31 \macc@nested@a\relax111{#1}%
  \else
  \def\gobble@till@marker##1\endmarker{}%
  \futurelet\first@char\gobble@till@marker#1\endmarker
  \ifcat\noexpand\first@char A\else \def\first@char{}%
  \fi
  \macc@nested@a\relax111{\first@char}%
  \fi
  \endgroup
}
\makeatother

%% file: abstract.tex


The Canonical Polyadic (CP) tensor decomposition is a well-known
method for interpretable analysis of high-dimensional
data~\cite{Kolda09TensorDecompositionsApplications}.  Recently, the
Generalized CP method ({\gcp}) was introduced by Hong and
Kolda~\cite{Hong20GeneralizedCanonicalPolyadic} to allow for flexible
choice of the loss function in the optimization problem defining the
CP model, enabling more interpretable decompositions of strongly
non-Gaussian data such as count or binary data.  Furthermore, Kolda
and Hong introduced in~\cite{Kolda20StochasticGradientsLargeScale} a
version of {\gcp} that leverages randomization and stochastic
optimization  to address scalability to large, sparse data sets.  In
this work, we take these ideas a step further and consider
synchronous and asynchronous algorithms for parallel GCP tensor
decomposition through the GenTen software package, exploiting both
shared and distributed memory parallelism.  We build on the shared
memory parallel CP decomposition algorithms utilizing Kokkos for
portability across CPU and GPU architectures described
in~\cite{Phipps19SoftwareSparseTensora} to support the random
sampling and stochastic optimization methods required by GCP.  We
then couple this approach to the well-known medium-grained
distributed memory parallelism scheme~\cite{SmithMediumGrained}
developed for traditional CP decompositions through MPI, providing a
synchronous, hybrid MPI+Kokkos, parallel GCP decomposition
capability.  Finally, we propose an asynchronous distributed
parallelism approach building on related techniques for federated
learning to achieve even better scalability to large data sets.  We
study the effectiveness of the proposed synchronous and asynchronous
approaches vis-\`{a}-vis computational cost and accuracy on synthetic
and publicly-available real-world datasets of varying sizes,
dimensions, and sparsity patterns using several loss functions.

%% file: intro.tex

Multi-way arrays called \textit{tensors} are useful for preserving
the natural structure and multilinear relationships found in the
high-dimensional data that is often of interest in data science
applications. Tensor decompositions based on the \textit{Canonical
Polyadic (CP)} model are used widely for varied tasks such as
online social networks \cite{GuPaPa17}, anomaly detection
\cite{FaGa16a}, compression of neural nets
\cite{JaSeAn15a,NoPoOsVe15,CoShSh16}, and health data analytics
\cite{WaChGhDe15},
especially when the data and corresponding model are
assumed to be low-rank. The Generalized CP
tensor decomposition ({\gcp}) introduced by Hong and
Kolda~\cite{Hong20GeneralizedCanonicalPolyadic} is a flexible
method for low-rank CP that allows for any arbitrary loss function using a
gradient-based optimization scheme. For large-scale tensor
decompositions in general and for {\gcp} in particular, scalability
remains a primary concern. Kolda and Hong introduced
in~\cite{Kolda20StochasticGradientsLargeScale} a version of {\gcp}
that leverages the {\adam} stochastic gradient
optimization method~\cite{Kingma15AdamMethodStochastic} to perform loss
function and gradient computations on $s$ samples of the input data
tensor so that the search path is computed from estimates of these
values, reducing the overall complexity to $\BigO(s)$.
Sequential implementations of {\gcp}-{\adam} are available in the
MATLAB and Python Tensor
Toolboxes~\cite{Bader17MATLABTensorToolbox,Dunlavy25PyttbPythonTensor}
that are suitable for small to medium sized data sets.

In this work, we build on these ideas to investigate parallel GCP
algorithms and associated software implementations.  Recognizing the
sequential {\gcp}-{\adam} approach described above employs similar
algorithmic building blocks as other traditional CP decompositions
methods, we use the shared-memory parallel CP decomposition
algorithms implemented within
GenTen\footnote{\href{https://www.github.com/sandialabs/genten}{https://www.github.com/sandialabs/genten}.}
as the foundation for constructing parallel GCP capabilities.  As
described in~\cite{Phipps19SoftwareSparseTensora}, GenTen leverages
Kokkos~\cite{Kokkos} to provide portability across contemporary CPU
and GPU computing architectures from a single algorithmic
implementation, yet still obtains performance that is commensurate
with tuned implementations for each architecture.  Here we leverage
similar techniques with Kokkos to construct shared-memory parallel
implementations of the additional computational kernels that are
unique to {\gcp}, providing a parallel {\gcp}-{\adam} algorithmic
implementation for CPU and GPU architectures within GenTen.
Furthermore, we extend this approach to synchronous distributed
parallelism with MPI to support {\gcp} decomposition of very large
tensor data sets with GenTen by incorporating the medium-grained
partitioning method~\cite{SmithMediumGrained} proposed by Smith and
Karypis for distributing sparse tensors across MPI processes.
However, due to the randomized nature of the {\gcp} algorithm, the
best partitioning method for the associated factor matrices is
unclear, so we investigate two approaches that give rise to very
different parallel communication patterns during the stochastic
gradient computation that exhibit differing cost and scalability properties.

Recently, increased attention has been paid to distributed settings
for GCP where local data sources are increasingly common (e.g.,
streaming sensor networks) or where privacy policies may forbid data
movement entirely (e.g., medical records). In these cases,
\textit{federated learning} has been employed to mitigate expensive collective
communications (e.g., \texttt{MPI\_Allreduce}) through asynchrony.
Asynchronous approaches leverage redundancy in computations and
data instances locally
to make progress and use collectives infrequently. Local
contributions are then combined through a reconciliation scheme to
form the global model.
Lewis and Phipps introduced
asynchronous {\gcp}-{\adam}
in~\cite{Lewis21LowCommunicationAsynchronousDistributed} utilizing
\textsc{ElasticAveraging}~\cite{Zhang15DeepLearningElastic} and
\textsc{LocalSGD}~\cite{Stich19LocalSGDConverges} for the
client-server coordination.
However, that work only demonstrated improved computational
performance and did not
consider the impact of asynchrony on convergence.
In this work, we extend a different federated technique, {\fedopt} described
in~\cite[Alg.
2]{Reddi21AdaptiveFederatedOptimization}, to asynchronous
{\gcp}-{\adam} with better convergence properties. In general,
federated learning algorithms pose
several challenges:
\begin{itemize}
  \item The communication frequency-update trade-off is complex. %
  \item Assumptions that guarantee an optimal convergence rate are
    difficult to enforce on real data. In our case, this is more
    complicated since the GCP optimization problem is nonconvex in general.
  \item Data heterogeneity across clients may lead to undesired
    convergence behavior. Federated learning algorithms may converge
    more quickly to solutions that overstate the importance of
    ``easier'' subsets of data when they reside on clients in
    isolation apart from ``harder'' subsets.
  \item The pseudo-gradients at the server level, i.e., estimates of
    the estimated gradients, are prone to high bias and variance
    which may impact the search path.
\end{itemize}
We consider two mitigations to address these challenges. First,
{\fedopt} more intelligently merges local gradients into the global
model than in prior work. In contrast to prior work, we study the impact of
asynchrony on convergence in addition to performance. Second, careful
parameter tuning has been demonstrated to account for bias and
variance in pseudo-gradients and the impact of data
heterogeneity, which was explored empirically and theoretically in
the original {\fedopt} paper from Reddi et al. Here we
advocate for an alternative to costly exhaustive grid search in
parameter tuning that is common in the uncertainty quantification
community yet underappreciated in the machine learning community,
namely Latin Hypercube
Sampling (LHS)~\cite{Helton03LatinHypercubeSampling}.

The contributions of this work are as follows:
\begin{itemize}
  \item We describe a first-ever parallel {\gcp} algorithm
    implemented within GenTen that supports hybrid shared-distributed
    memory parallelism through MPI+Kokkos, including demonstrations
    on a multi-node cluster including both CPU and GPU architectures.
  \item We propose a fused sampling approach that combines the
    sampled tensor and sampled gradient computations within {\gcp}
    that reduces computational cost by reducing memory traffic.
  \item We propose a novel federated learning approach for GCP that
    incorporates asynchronous parallelism to reduce synchronization
    costs and facilitate local computation, resulting in reduced
    computational costs, increased scalability, and increased privacy.
\end{itemize}
Background on {\gcp} and related work on parallel tensor
decomposition is provided in~\cref{background}.  Then parallelization
of {\gcp} using Kokkos and MPI is described
in~\cref{genten-dist-par}, including the fused sampling approach in
\cref{sec:stoch_gradient}.  Next, the proposed asynchronous federated
learning approach is described in~\cref{fed-opt}, followed by several
computational experiments demonstrating performance, portability, and
scalability of the synchronous and asynchronous approaches on CPU and
GPU architectures in~\cref{experiments}.  Finally, conclusions and
future work are summarized in~\cref{conc}.


%% file: background.tex
\subsection{Notation}
A tensor $\tns{X}$ is a $d$-way array with dimensions $I_1 \times I_2
\times \dots \times I_d$.
Following standard practice, we denote tensors by bold calligraphic
letters (e.g., $\tns{X}$),
matrices by bold capital letters ($\mat{A}$),
vectors by bold lowercase letters ($\mat{a}$)
and scalars by lowercase letters ($a$).
We use multi-index notation to indicate tensor elements, i.e.,
$x_i \equiv x_{i_1\dots i_d}$ denotes the entry $i=(i_1,\dots,i_d)$
of $\tns{X}$.
A tensor
\textit{slice} is a ($d-1$)-dimensional subtensor given by fixing one
index, i.e., $\tns{X}_{i_1,:,\dots,:}$. A
tensor \textit{fiber} is a one-dimensional subtensor analogous to a
column or row vector given by fixing all but one index, e.g.,
$\mat{x}_{:,i_2,\dots,i_d}$.

The
\emph{rank-$R$ canonical polyadic (CP) tensor model of
$\tns{X}$}~\cite{Harshman70FoundationsPARAFACProcedure,
Carroll70AnalysisIndividualDifferences} is:
\begin{equation}
  \label{eq:cp}
  \tns{X} \approx \tns{M}
  = \llbracket \mat{A}^{(1)}, \ldots, \mat{A}^{(d)} \rrbracket
  \coloneqq \sum_{r=1}^R \mat{a}^{(1)}_r \circ \dots
  \circ \mat{a}^{(d)}_r,
\end{equation}
where each $\mat{A}^{(k)} \in \R^{I_k \times R}$ is a \emph{factor
matrix} with $I_k$ rows
and $R$ columns, the $j$-th \emph{component} of the mode-$k$ factor matrix is
the column vector $\mat{a}^{(k)}_j$, and $\circ$ denotes the
$d$-way outer product of $d$ factors.
We refer to the form $\tns{M} = \llbracket\mat{A}^{(1)}, \ldots,
\mat{A}^{(d)} \rrbracket$ as a \emph{Kruskal tensor}.

\subsection{Generalized Canonical Polyadic Decomposition (GCP)}

The Generalized Canonical Polyadic (GCP) low-rank tensor
decomposition, introduced by Hong and
Kolda~\cite{Hong20GeneralizedCanonicalPolyadic}, enables a more
flexible choice of loss
functions beyond traditional sum-of-squared error for data where a
Gaussian statistical model is not appropriate by solving a nonlinear,
nonconvex maximum likelihood
estimation problem:
\begin{equation} \label{eq:gcp-model}
  \min_{\tns{M}} F(\tns{X}, \tns{M}) = \sum_{\mind} f (x_{\mind}, m_{\mind}).
\end{equation}
GCP takes advantage of the relationship between the mode-$k$
unfolding of a Kruskal tensor, $\mat{M}_{(k)}$, and the Khatri-Rao
product $\mat{Z}_k$ of its factor matrices. This fact is used to
define the GCP gradient for mode-$k$ as a matricized tensor times
Khatri-Rao product (MTTKRP):
\begin{equation} \label{eq:gcp-gradient}
  \mat{G}^{(k)} = \frac{\partial F}{\partial \mat{A}^{(k)}} =
  \tns{Y}_{(k)} \mat{Z}_k^T
\end{equation}
where $ \mat{Z}_k \equiv \mat{A}^{(d)} \odot \cdots \odot \mat{A}^{(k+1)}
\odot \mat{A}^{(k)} \odot \cdots \odot \mat{A}^{(1)}$ and $y_{\mind} =
\frac{\partial f}{\partial m}(x_{\mind}, m_{\mind})$. MTTKRP is a
fundamental computational kernel in many tensor decompositions and
its use in~\cref{eq:gcp-gradient} motivates several first-order
methods to solve~\cref{eq:gcp-model}, including quasi-Newton methods
such as Limited-Memory BFGS and stochastic
gradient descent ({\sgd}).

\subsection{ADAM for Stochastic Gradient Descent}
{\adam}~\cite{Kingma15AdamMethodStochastic} is a popular variant of
{\sgd} that adjusts the learning rate dynamically and uses momentum
to escape local minima that are far from the maximum likelihood
estimator. For reference, we provide the model update performed in each iteration of
{\adam} (adapted from~\cite{Hong20GeneralizedCanonicalPolyadic}) where the model and gradient are Kruskal tensors
 (here $\ast$ and $\oslash$ denote
element-wise multiplication and division, respectively).

\begin{algorithm}
  \caption{\textsc{Adam} update}
  \label{alg:adam-update}
  \begin{algorithmic}[1]
    \Function{update}{Model $\tns{M} = \llbracket \mat{A}^{(1)}, \ldots,
      \mat{A}^{(d)} \rrbracket$, stochastic gradient $\tns{U} = \llbracket
    \mat{G}^{(1)}, \ldots, \mat{G}^{(d)} \rrbracket$}
    \State \textbf{global variables} $\{\mat{B}^{(k)}\}, \{\mat{C}^{(k)}\},
    \; k = 1, \ldots, d$.
    \For{$k = 1, \ldots, d$}
    \State $\mat{B}^{(k)} \gets \beta_1 \mat{B}^{(k)} +
    (1-\beta_1)\mat{G}^{(k)}$
    \State $\mat{C}^{(k)} \gets \beta_2 \mat{C}^{(k)} +
    (1-\beta_2)\mat{G}^{(k)}\ast\mat{G}^{(k)}$
    \State $\mat{\hat{B}}^{(k)} \gets \mat{B}^{(k)} / (1-\beta_1^t)$
    \State $\mat{\hat{C}}^{(k)} \gets \mat{C}^{(k)} / (1-\beta_2^t)$
    \State $\mat{A}^{(k)} \gets \mat{A}^{(k)} - \alpha (\mat{\hat{B}}^{(k)}
    \oslash \sqrt{\mat{\hat{C}}^{(k)} + \epsilon})$
    \State $\mat{A}^{(k)} \gets \max\{ \mat{A}^{(k)}, l \}$
    \Comment{$l= \text{lower bound}$}
    \EndFor
    \EndFunction
  \end{algorithmic}
\end{algorithm}%

\subsection{Related Work}
{\gcp} was introduced by Hong and Kolda for dense tensors
in~\cite{Hong20GeneralizedCanonicalPolyadic} using quasi-Newton
optimization while an approach based on Gauss-Newton optimization was
considered in~\cite{VaVeLa21}.  {\gcp} was extended to sparse tensors
through stochastic optimization methods
in~\cite{Kolda20StochasticGradientsLargeScale}.  An alternative
approach based on fiber sampling and stochastic mirror descent was
described in~\cite{PuIbFuHo22}.  For Poisson loss, a hybrid
stochastic-deterministic algorithm was studied
in~\cite{Myers25TensorDecompositionsCount}.  {\gcp} was extended to
the streaming setting
in~\cite{Phipps23StreamingGeneralizedCanonicala} and included a
Kokkos-based shared-memory parallel implementation in GenTen.

Beyond the streaming work cited above, we believe this is the first
work to consider either shared- or distributed-memory parallel {\gcp}
for sparse tensors.  
Conversely, numerous works have considered parallel CP decomposition
using the traditional sum-of-squares loss function for Gaussian data.
Most of this work focuses on sparse tensor storage formats to reduce
costs for the Matricized Tensor Times Khatri-Rao Product (MTTKRP)
kernel that is critical for performance, such as compressed sparse
fiber~\cite{SmRaSiKa15}, dimension trees~\cite{Kaya:2018dw}, flagged
coordinate~\cite{2017arXiv170509905L}, block
coordinate~\cite{Li:2018sc} and linearized indices~\cite{alto2021}.
Here we focus on the traditional sparse coordinate format given its
mode-agnostic storage, potential for high performance as demonstrated
in~\cite{Phipps19SoftwareSparseTensora}, and low setup costs for
sparse sampling.  For distributed parallelism of traditional CP
decomposition, three primary approaches have emerged for distributing
the tensor and factor matrices across processors.  Tools such as
DFacTo~\cite{ChVi14} implement a coarse-grained approach given by a
one-dimensional distribution of tensor slices across processors while
factor matrices are replicated across processors.  This approach
limits parallelism to a single mode and incurs high storage costs to
the replicated factor matrices.  Conversely,
HyperTensor~\cite{KaUc15b} proposes a fine-grained distribution of
tensor nonzeros across processors using hypergraph partitioning.
While this approach results in optimal communication volume, it has
substantial setup costs to compute the hypergraph partitioning.
Instead, we follow the medium-grained partitioning approach presented
in~\cite{SmithMediumGrained} which distributes tensors across
processors using a $d$-way processor grid since it does not incur
expensive preprocessing but still provides substantial parallelism.
How to distribute the factor matrices in a manner suitable for {\gcp}
is studied in this paper.  Finally,
BalaCPD~\cite{Miao22BALACPDBALancedAsynchronous}  considers load
balancing using the medium-grained approach and also proposes a
one-sided communication scheme within MTTKRP using OpenSHMEM.

The development of federated learning techniques for {\gcp} has roots
in distributing stochastic gradient descent ({\sgd}) computations.
Dean et al.~\cite{Dean12LargeScaleDistributed} introduced
\textsc{DownpourSGD} to leverage asynchrony which reduces costly
collective communications by exploiting locality. The trade-off is
that this incurs overhead due to redundant computations and replicated
data on each shard in the processor topology.

Subsequent works have proposed various synchronization schemes that
are implemented in GenTen.
\textsc{LocalSGD}~\cite{Stich19LocalSGDConverges} uses simple
averaging to incorporate local models globally. \textit{Elastic
Averaging (EA)}~\cite{Zhang15DeepLearningElastic} updates the global
model with local moving averages. Both schemes tailored for GCP
appeared in~\cite{Lewis21LowCommunicationAsynchronousDistributed} and
are implemented in GenTen. Synchronous {\gcp}-{\sgd} and
\textsc{LocalSGD} are presented in \Cref{alg:gcp-sgd-iter} and
\Cref{alg:local-sgd-iter} for reference, respectively.
Note that Line~\ref{alg:gcp-sgd-iter:gradient}
in~\Cref{alg:gcp-sgd-iter}  and
Line~\ref{alg:local-sgd-iter:gradient} in~\Cref{alg:local-sgd-iter}
are implementation details in GenTen that applies the MTTKRP and
computes the gradient factor matrices in~\cref{eq:gcp-gradient}. They
are given more extensive treatment in~\cref{sec:gcp_kokkos}.
These
asynchronous approaches have several limitations. Reddi et
al.~\cite{Reddi21AdaptiveFederatedOptimization} cite convergence
issues with \textsc{LocalSGD} due to client drift and a lack of
adaptivity in the server optimization, which is necessary for
heterogeneous data. Those authors address adaptivity and introduce
{\fedopt}, a generalized framework for federated learning that
permits {\adam}-type updates at both the client and server levels.

\begin{algorithm}
  \caption{\textsc{Adam} mini-batch iteration (synchronous)}
  \label{alg:gcp-sgd-iter}
  \begin{algorithmic}[1]
    \Function{AdamSGDIteration}{$\tns{M}$, $sampler$}
    \State $\tns{G} \gets sampler.Gradient(\tns{M})$
    \label{alg:gcp-sgd-iter:gradient}
    \State $Allreduce(\tns{G})$ \label{alg:gcp-sgd-iter:allreduce}
    \State $\tns{M}.update(\tns{G})$
    \EndFunction
  \end{algorithmic}
\end{algorithm}%

\begin{algorithm}
  \caption{Local SGD mini-batch iteration (asynchronous)}
  \label{alg:local-sgd-iter}
  \begin{algorithmic}[1]
    \Function{LocalSGDIteration}{$\tns{M}$, $sampler$, epoch $e > 0$,
    asynchrony $\tau > 0$}
    \If{$\text{mod}(e, \tau) = 0$}
    \State $Allreduce(\tns{M})$
    \State $\tns{M} \gets \tns{M} / numMPIRanks$
    \EndIf
    \State $\tns{G} \gets sampler.Gradient(\tns{M})$
    \label{alg:local-sgd-iter:gradient}
    \State $\tns{M}.update(\tns{G})$
    \EndFunction
  \end{algorithmic}
\end{algorithm}%

%% file: genten-distributed-parallelism.tex

The GCP method introduced
in~\cite{Hong20GeneralizedCanonicalPolyadic} and extended to sparse
tensors in~\cite{Kolda20StochasticGradientsLargeScale} focused on a
serial implementation in the Matlab Tensor
Toolbox~\cite{Bader17MATLABTensorToolbox}.  Here we described a
hybrid-parallel implementation within
GenTen~\cite{Phipps19SoftwareSparseTensora} that mirrors the Tensor
Toolbox implementation but leverages Kokkos~\cite{Kokkos} for
performance portable shared memory parallelism and MPI for
distributed parallelism.  We first describe our approach to threaded
parallelism in GCP using Kokkos in \cref{sec:gcp_kokkos} and then
extend that implementation with distributed parallelism using MPI in
\cref{sec:gcp_mpi}.

\subsection{Shared Memory Parallelism with Kokkos}
\label{sec:gcp_kokkos}

Threaded parallelism in GenTen leveraging Kokkos for CP decomposition
of sparse tensors was extensively described
in~\cite{Phipps19SoftwareSparseTensora}.  Kokkos is a C++ library
that provides data structures and interfaces that applications can
use to express fine-grained threaded parallelism in a manner that is
portable across contemporary CPU and GPU architectures.  That work
was relevant to the CP-ALS algorithm suitable for least-squares loss
and focused on leveraging Kokkos for parallelizing the MTTKRP
operation required to update factor matrix components at each
iteration.  Applying similar techniques to parallelize the GCP
algorithm reviewed above with Kokkos requires parallel
implementations of three high-level phases of the calculation:
stochastic sampling to construct the gradient tensor $\tns{Y}$,
applying MTTKRP to compute the gradient factor matrices
$\{\mat{G}_k\}$, and updating the factor matrices using {\adam}.
Each of these is discussed below.

\subsubsection{Parallel Stochastic Sampling}
\label{sec:stoch_sampling}

The tensor $\tns{Y}$ appearing in \cref{eq:gcp-gradient} (heretofore
called the gradient tensor) is in general dense even if the original
data tensor $\tns{X}$ is sparse, limiting scalability to large sparse
tensors.  To overcome this challenge, the authors
in~\cite{Kolda20StochasticGradientsLargeScale} turned to stochastic
gradient descent as a solution strategy, and proposed several
sampling algorithms that produce a stochastic approximation
$\tilde{\tns{Y}}\approx\tns{Y}$.  The gradient tensor $\tns{Y}$ in
\cref{eq:gcp-gradient} is then replaced by $\tilde{\tns{Y}}$
resulting in the stochastic gradient $\tilde{\mat{G}}^{(k)} =
\tilde{\tns{Y}}_{(k)} \mat{Z}_k^T$.  Several sampling algorithms were
introduced in~\cite{Kolda20StochasticGradientsLargeScale}, but here
we focus on the two most relevant for sparse tensors:  stratified and
semi-stratified sampling.

Stratified sampling recognizes that the probability of selecting
nonzeros becomes vanishingly small as the tensor becomes larger and
sparser, and therefore proposes that zeros and nonzeros be sampled
separately.  Following~\cite{Kolda20StochasticGradientsLargeScale},
given non-negative integers $p$ and $q$, we sample {\em with
replacement} $p$ nonzeros of $\tns{X}$ and $q$ zeros.  GenTen stores
sparse tensors using the coordinate format consisting of a list of
$N$ tuples $(i_1,\dots,i_d)$ storing the coordinates of each nonzero
and a list of $N$ floating-point values $x_{i_1,\dots,i_d}$ providing
tensor entries.  Thus sampling zeros is straightforward and just
requires randomly selecting integers in the range $[1,\dots,N]$ and
reading the corresponding coordinates and value.  A new nonzero is
then added to $\tilde{\tns{Y}}$ defined by $\tilde{y}_{i_1,\dots,i_d}
= \frac{N}{p}\frac{\partial f}{\partial
m}(x_{i_1,\dots,i_d},m_{i_1,\dots,i_d})$.  Sampling zeros is more
challenging, however, because zeros are not explicitly stored.
Instead, for each zero sample, $d$ integers $i_1,\dots,i_d$ are each
randomly chosen out of the ranges
$[1,\dots,I_1],\dots,[1,\dots,I_d]$.  The tensor is then searched to
determine if $(i_1,\dots,i_d)$ is contained within the list of
nonzeros.  If it is not found, the sample corresponds to a zero,
otherwise the process repeats until the sampled entry is not found.
A new nonzero is then added to $\tilde{\tns{Y}}$ given by
$\tilde{y}_{i_1,\dots,i_d} = \frac{M-N}{q}\frac{\partial f}{\partial
m}(0,m_{i_1,\dots,i_d})$, where $M \equiv \prod_{j=1}^d I_j$.

Stratified sampling is implemented in GenTen by parallelizing the
loops over nonzero and zero samples using the Kokkos parallel-for
pattern.  Each thread computes entries in $\tilde{\tns{Y}}$ as
described above for both nonzeros and zeros, which requires computing
CP model entries $m_{i_1,\dots,i_d} = \sum_{j=1}^R
a^{(1)}_{i_1,j}\dots a^{(d)}_{i_d,j}$.  Since the number of nonzeros
in the stochastic gradient tensor $\tilde{\tns{Y}}$ is given by
$p+q$, space for storing the coordinates and values is preallocated,
and thus each thread just needs to assign entries in the list of
nonzeros without synchronization.\footnote{Because samples are
  generated with replacement, different threads, or even the same
  thread, might sample a zero/nonzero more than once, in which case
  that entry appears multiple times in the $\tilde{\tns{Y}}$ tensor.
However, those entries will have the same value by construction.}  To
speedup the required search for each candidate zero, GenTen sorts the
entries lexicographically based on their coordinate index, which
allows each thread to search the tensor in $O(\log N)$ operations.
GenTen also provides a hash map capability to store tensor nonzeros
using \texttt{Kokkos::unordered\_map} using nonzero coordinate
indices as keys, which increases setup costs but then allows each
thread to search the tensor in $O(1)$ operations.

Semi-stratified sampling was introduced
in~\cite{Kolda20StochasticGradientsLargeScale} as an improvement to
stratified sampling that eliminates the need to search the tensor for
each candidate zero sample.  The idea is to still samples nonzeros
and zeros separately, but assume each zero sample corresponds to a
zero in the tensor and then correct the nonzero samples for the case
when a nonzero was incorrectly sampled as a zero.  The algorithmic
implementation is the same as described above, without the search,
but with stochastic gradient tensor entries for nonzero samples given
by $\tilde{y}_{i_1,\dots,i_d} = \frac{N}{p}\left(\frac{\partial
  f}{\partial m}(x_{i_1,\dots,i_d},m_{i_1,\dots,i_d}) - \frac{\partial
f}{\partial m}(0,m_{i_1,\dots,i_d})\right)$ (the entries for zero
samples are unchanged).

\subsubsection{Parallel Stochastic Gradient}
\label{sec:stoch_gradient}

Once the gradient tensor $\tilde{\tns{Y}}$ has been computed, the
stochastic gradient $\tilde{\mat{G}}^{(k)} = \tilde{\tns{Y}}_{(k)}
\mat{Z}_k^T$ can then be computed using the Kokkos-based MTTKRP
kernels as described in~\cite{Phipps19SoftwareSparseTensora}.  Note,
however, that a new gradient tensor $\tilde{\tns{Y}}$ is computed for
each stochastic gradient descent iteration, each of which has a
potentially different sparsity pattern.  Therefore, the
permutation-based MTTKRP that was found to be the most effective for
CP-ALS in~\cite{Phipps19SoftwareSparseTensora} is usually not the
most effective for GCP since it requires computing new permutation
arrays for each sampled tensor, and the cost of computing those
permutation arrays is comparable to the MTTKRP itself.  Thus, the
most effective MTTKRP algorithm for GCP is architecture and problem
dependent.  For recent GPUs with hardware atomic support, the
atomic-based algorithm described
in~\cite{Phipps19SoftwareSparseTensora} is the best choice.  For
CPUs, either the atomic approach or an approach based on thread
privatization where each thread updates a thread-local copy of the
factor matrices can be more efficient, depending on the problem size
and shape, with atomics often being more efficient for tensors with
large mode sizes and thread-privatization otherwise.  GenTen also
provides a serial approach for use in serial or MPI-only contexts
that eliminates any need for managing thread race conditions (note
  that because of this, MPI+Serial can often be more efficient than
OpenMP-only on CPUs).

It is worth noting, however, the structure of the sampling operations
described above is quite similar to MTTKRP applied to
$\tilde{\tns{Y}}$ as both involve looping over all entries in
$\tilde{\tns{Y}}$ and accessing the same factor matrix entries.
Thus, performance can be potentially improved by fusing these two
operations, whereby for each sampled nonzero/zero, the contribution
of that sampled entry to the MTTKRP calculation $y_i
a^{(1)}_{i_1,j}\dots a^{(d)}_{i_d,j}$ is immediately computed,
without explicitly constructing the gradient tensor
$\tilde{\tns{Y}}$.  This approach has several advantages:  it
eliminates the additional memory cost of storing $\tilde{\tns{Y}}$,
saves the memory bandwidth costs of writing it to main memory during
sampling and then reading it in MTTKRP (since it is often too large
to fit in cache), and reduces memory bandwidth costs of accessing the
factor matrix coefficients by only reading them from main memory once
instead of twice.  This is referred to as the Fused Sampling-MTTKRP
algorithm, and is supported for all Kokkos backends using atomics,
thread privatization, or serial, but only for the semi-stratified
sampling approach.

\subsubsection{Parallel {\adam}}
\label{sec:parallel_adam}

Kokkos parallelization of the {\adam} update algorithm described in
\Cref{alg:adam-update} is straightforward whereby factor matrix
operations in each step of that algorithm are mapped across threads
using the Kokkos parallel-for pattern.  Naively, this would require a
separate parallel-for kernel for each tensor dimension, e.g., a
separate invocation of  a parallel kernel implementing $\mat{B}^{(k)}
\gets \beta_1 \mat{B}^{(k)} + (1-\beta_1)\mat{G}^{(k)}$ for each $k$.
However, GenTen stores the coefficients for all factor matrices in a
single contiguous array for the GCP algorithm, and provides views to
subsets of that array for each factor matrix.  This allows all factor
matrices to be updated with a single kernel invocation, e.g.,
$\mat{b} \gets \beta_1 \mat{b} + (1-\beta_1)\mat{g}$, and facilitates
a simpler implementation as it only operates on vectors instead of
collections of matrices.

\subsection{Distributed Parallelism with MPI}
\label{sec:gcp_mpi}

In addition to shared memory parallelism with Kokkos, GenTen also
supports distributed memory parallelism using MPI where the tensor
and factor matrices are distributed across MPI processes.  These two
modes of parallelism can be combined to form hybrid parallelism where
each MPI process executes the GCP operations described above in
parallel using Kokkos.  On CPUs, this typically amounts to binding
MPI processes to a collection of cores and using Kokkos with the
OpenMP backend for thread parallelism amongst those cores.  For a
multicore node with $n$ cores, virtually any combination of MPI
processes and threads is typically supported, but the most performant
combinations typically are one MPI process per node and $n$ OpenMP
threads, one MPI process per CPU socket (or more generally, NUMA
domain) and $n/s$ threads where $s$ is the number of sockets, and $n$
MPI processes with no threading (typically using the Kokkos Serial
backend), each bound to a core.  In all cases, OpenMP threads should
be bound to their respective cores to make best use of the cache
hierarchy.  While Kokkos and GenTen support the use of hyperthreads,
it typically not advisable since it most often results in reduced
performance.  For GPU architectures, the recommended approach is to
use one MPI process per GPU on each node.  GenTen assumes the chosen
MPI library is ``GPU-aware'' and supports communicating data directly
from/to GPU memory buffers without explicitly transferring data
to/from the host CPU memory space first.

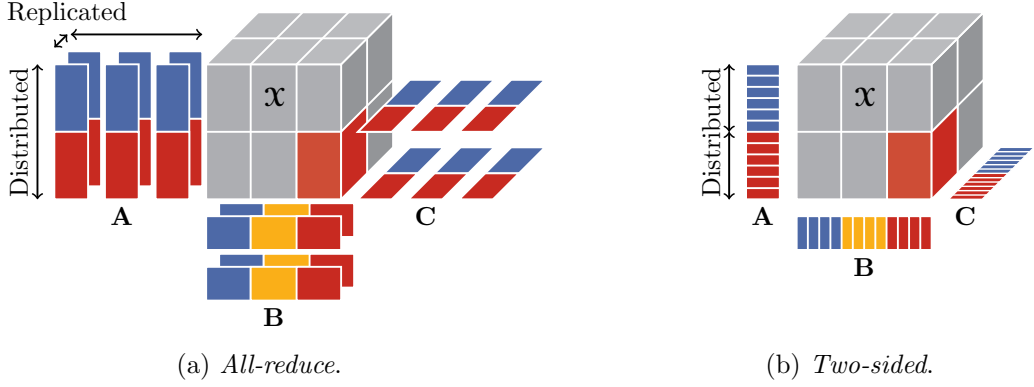
\begin{figure}
  \centering
  \begin{subfigure}[t]{0.5\textwidth}
    \centering
    \begin{adjustbox}{width=\linewidth}
      \input{images/genten_allreduce}
    \end{adjustbox}
    \caption{\textit{All-reduce}.}
    \label{fig:genten-distributed-parallelism:allreduce}
  \end{subfigure}\hfill
  \begin{subfigure}[t]{0.5\textwidth}
    \centering
    \begin{adjustbox}{width=\linewidth}
      \input{images/genten_two_sided}
    \end{adjustbox}
    \caption{\textit{Two-sided}.}
    \label{fig:genten-distributed-parallelism:tpetra}
  \end{subfigure}
  \caption{Two versions of distributed parallelism in GenTen. The
    figures illustrate a tensor $\tns{X}$ and Kruskal tensor $\tns{M} =
    \llbracket \mat{A}, \mat{B}, \mat{C} \rrbracket$ partitioned across
    12 processor ranks. \textit{Distributed} factor matrices are those
    associated with a given mode $k$ of the tensor partitioned across
    the associated mode-$k$ fiber communicators. \textit{Replicated}
    factor matrices are those distributed across the associated slice
  communicators.}
  \label{fig:genten-distributed-parallelism}
\end{figure}

For distributing the tensor in parallel across MPI processes, GenTen
follows the medium-grained partitioning approach described
in~\cite{SmithMediumGrained} whereby each dimension $k$ of the tensor
is partitioned among $N_k$ processes where $N \equiv \prod_{j=1}^d
N_j $ is the total number of MPI processes.  This results in a
partitioning of the tensor into blocks according to a $d$-way
processor grid $N_1\times\dots\times N_d$, as demonstrated in
\Cref{fig:genten-distributed-parallelism}.  For a user-chosen number
of processes $N$, GenTen computes the per-dimension processor counts
$N_k$ to minimize the total storage of the factor matrices by
considering all possible combinations of divisors of $N$, assuming
the ``all-reduce'' factor matrix distribution/MTTKRP approach
described below.  Thus, the MPI communicator is a Cartesian
communicator which supports efficient communication between
rectilinear subsets of processors.  In particular, we define slice
and fiber sub-communicators analogously to slices and fibers in the
tensor, i.e., given a processor with processor multi-index
$(i_1,i_2,\dots,i_d)$, the mode-1 slice and fiber sub-communicators
for this processor consist of all processors with multi-index
$(i_1,:,\dots,:)$ and $(:,i_2,\dots,i_d)$, respectively.

Given the geometric partitioning of the tensor described above,
GenTen supports two approaches for distributing factor matrices in a
conformal manner.  The first approach, which we call the
``all-reduce'' approach due to its implied parallel implementation of
distributed MTTKRP, distributes factor matrices associated with a
given mode $k$ of the tensor across the associated mode-$k$ fiber
communicators and replicates those factor matrices across the
associated slice communicators.  An example for a 3-way tensor is
shown in \Cref{fig:genten-distributed-parallelism:allreduce}.  This
is the approach that was originally implemented in GenTen and the
processor grid factorization described above minimizes the total
factor matrix storage.  An advantage of this approach is each
processor has access to
the data it requires for its local computations. As a result,
computing the MTTKRP is relatively straightforward. When the global
model is updated, the collective \texttt{MPI\_Allreduce}, is
called across the corresponding slice communicator. For {\sgd} in general and
{\adam} in particular, GCP leverages this distributed parallelism
model to sample tensor and gradient entries local only to a given
processor. However, this approach has two disadvantages.  First, the
collective \texttt{MPI\_Allreduce} can be expensive since the communication
volume scales with tensor dimensions, i.e., the update is dense
regardless of sparsity pattern.  Second, it only provides limited
parallelism in factor matrix computations such as the {\adam} updates.

The second approach supported by GenTen is referred to as the
``two-sided'' approach (owing to the two-sided communication pattern
described next), which was originally described
in~\cite{SmithMediumGrained} and implemented to overcome these
disadvantages.  It distributes factor matrices across all processors,
but is done in a manner conformal to the partitioning of the tensor.
For a factor matrix associated with a given mode $k$, the factor
matrix is first partitioned across $N_k$ blocks of processors where
the processor boundaries are determined by the corresponding
boundaries for the tensor.  Then each block is further partitioned
uniformly across the remaining $N/N_k$ processors.  In this way, each
processor still only needs to communicate with other processors
contained within its slice sub-communicator for the parallel MTTKRP
computation.  However, each processor must receive off-processor
factor matrix entries before executing the local MTTKRP computation
and then send its contributions for rows of the result factor matrix
it does not own to the owning processors (both of which can be
  implemented efficiently by \texttt{MPI\_Alltoall} for the
corresponding slice sub-communicator).  This approach increases
factor matrix parallelism and reduces communication volume in MTTKRP
compared to the all-reduce approach since only the necessary factor
matrix rows are communicated between processors, but has the
disadvantage of requiring additional communication and overhead to
determine which rows will be communicated to which processors based
on the sparsity pattern of the tensor.  This is particularly
problematic for GCP since the sparsity pattern changes for each
sampled gradient tensor.  Thus which of these two approaches is more
efficient is problem dependent based on the size, shape, and sparsity
of the tensor as well as the number of chosen gradient samples.

Distributed parallel sampling and {\adam} updates are purely local
operations that do not require parallel communication and thus are
straightforward.  For stratified/semi-stratified sampling, the number
of zero and nonzero samples per processor are computed by uniformly
distributing the total number of zero and nonzero samples across all processors.

%% file: images/genten_allreduce.tex
\newcommand{\Depth}{2}
\newcommand{\Height}{2}
\newcommand{\Width}{2}
\newcommand{\Offset}{0.5}
\begin{tikzpicture}

  \coordinate (origin) at (0,0,0);
  \coordinate (lowerbackright) at (\Width,0,0);
  \coordinate (upperbackright) at (\Width,\Height,0);
  \coordinate (upperbackleft) at (0,\Height,0);

  \coordinate (lowerfrontleft) at (0,0,\Depth);
  \coordinate (lowerfrontright) at (\Width,0,\Depth);
  \coordinate (upperfrontright) at (\Width,\Height,\Depth);
  \coordinate (upperfrontleft) at (0,\Height,\Depth);

  \input{images/tensor3d.tex}

  
  \coordinate(rui0) at (\Width+0.5*\Offset,0+0.5*\Height,\Depth);
  \coordinate(rui1) at (\Width+0.5*\Offset,0+0.5*\Height,0);
  \coordinate(rui2) at (\Width+1.5*\Offset,0+0.5*\Height,0);
  \coordinate(rui3) at (\Width+1.5*\Offset,0+0.5*\Height,\Depth);
  \draw[white, thick, fill=myblue] (rui0) -- (rui1) -- (rui2) -- (rui3) -- cycle;
  \draw[white, thick] ($(rui0)!0.5!(rui1)$) -- ($(rui2)!0.5!(rui3)$);

  \coordinate(rli0) at (\Width+0.5*\Offset,0,\Depth);
  \coordinate(rli1) at (\Width+0.5*\Offset,0,0);
  \coordinate(rli2) at (\Width+1.5*\Offset,0,0);
  \coordinate(rli3) at (\Width+1.5*\Offset,0,\Depth);
  \draw[white, thick, fill=myblue] (rli0) -- (rli1) -- (rli2) -- (rli3) -- cycle;
  \draw[white, thick] ($(rli0)!0.5!(rli1)$) -- ($(rli2)!0.5!(rli3)$);

  \coordinate(rum0) at (\Width+2.0*\Offset,0+0.5*\Height,\Depth);
  \coordinate(rum1) at (\Width+2.0*\Offset,0+0.5*\Height,0);
  \coordinate(rum2) at (\Width+3.0*\Offset,0+0.5*\Height,0);
  \coordinate(rum3) at (\Width+3.0*\Offset,0+0.5*\Height,\Depth);
  \draw[white, thick, fill=myblue] (rum0) -- (rum1) -- (rum2) -- (rum3) -- cycle;
  \draw[white, thick] ($(rum0)!0.5!(rum1)$) -- ($(rum2)!0.5!(rum3)$);

  \coordinate(rlm0) at (\Width+2.0*\Offset,0,\Depth);
  \coordinate(rlm1) at (\Width+2.0*\Offset,0,0);
  \coordinate(rlm2) at (\Width+3.0*\Offset,0,0);
  \coordinate(rlm3) at (\Width+3.0*\Offset,0,\Depth);
  \draw[white, thick, fill=myblue] (rlm0) -- (rlm1) -- (rlm2) -- (rlm3) -- cycle;
  \draw[white, thick] ($(rlm0)!0.5!(rlm1)$) -- ($(rlm2)!0.5!(rlm3)$);

  \coordinate(ruo0) at (\Width+3.5*\Offset,0+0.5*\Height,\Depth);
  \coordinate(ruo1) at (\Width+3.5*\Offset,0+0.5*\Height,0);
  \coordinate(ruo2) at (\Width+4.5*\Offset,0+0.5*\Height,0);
  \coordinate(ruo3) at (\Width+4.5*\Offset,0+0.5*\Height,\Depth);
  \draw[white, thick, fill=myblue] (ruo0) -- (ruo1) -- (ruo2) -- (ruo3) -- cycle;
  \draw[white, thick] ($(ruo0)!0.5!(ruo1)$) -- ($(ruo2)!0.5!(ruo3)$);

  \coordinate(rlo0) at (\Width+3.5*\Offset,0,\Depth);
  \coordinate(rlo1) at (\Width+3.5*\Offset,0,0);
  \coordinate(rlo2) at (\Width+4.5*\Offset,0,0);
  \coordinate(rlo3) at (\Width+4.5*\Offset,0,\Depth);
  \draw[white, thick, fill=myblue] (rlo0) -- (rlo1) -- (rlo2) -- (rlo3) -- cycle;
  \draw[white, thick] ($(rlo0)!0.5!(rlo1)$) -- ($(rlo2)!0.5!(rlo3)$);

  \draw[white, thick, fill=myred] (rlo0) -- ($(rlo0)!0.5!(rlo1)$) --
  ($(rlo2)!0.5!(rlo3)$) -- (rlo3) -- (rlo0) -- cycle;
  \draw[white, thick, fill=myred] (rlm0) -- ($(rlm0)!0.5!(rlm1)$) --
  ($(rlm2)!0.5!(rlm3)$) -- (rlm3) -- (rlm0) -- cycle;
  \draw[white, thick, fill=myred] (rli0) -- ($(rli0)!0.5!(rli1)$) --
  ($(rli2)!0.5!(rli3)$) -- (rli3) -- (rli0) -- cycle;
  \draw[white, thick, fill=myred] (ruo0) -- ($(ruo0)!0.5!(ruo1)$) --
  ($(ruo2)!0.5!(ruo3)$) -- (ruo3) -- (ruo0) -- cycle;
  \draw[white, thick, fill=myred] (rum0) -- ($(rum0)!0.5!(rum1)$) --
  ($(rum2)!0.5!(rum3)$) -- (rum3) -- (rum0) -- cycle;
  \draw[white, thick, fill=myred] (rui0) -- ($(rui0)!0.5!(rui1)$) --
  ($(rui2)!0.5!(rui3)$) -- (rui3) -- (rui0) -- cycle;

  \coordinate (lbi0) at (-0.5*\Offset,0,-1*\Offset + \Depth);
  \coordinate (lbi1) at (-1.5*\Offset,0,-1*\Offset + \Depth);
  \coordinate (lbi2) at (-1.5*\Offset,\Height,-1*\Offset + \Depth);
  \coordinate (lbi3) at (-0.5*\Offset,\Height,-1*\Offset + \Depth);
  \draw[white, thick, fill=myblue] (lbi0) -- (lbi1) -- (lbi2) -- (lbi3) -- cycle;
  \draw[white, thick] ($(lbi0)!0.5!(lbi3)$) -- ($(lbi1)!0.5!(lbi2)$);

  \coordinate (lbm0) at (-2.0*\Offset,0,-1*\Offset + \Depth);
  \coordinate (lbm1) at (-3.0*\Offset,0,-1*\Offset + \Depth);
  \coordinate (lbm2) at (-3.0*\Offset,\Height,-1*\Offset + \Depth);
  \coordinate (lbm3) at (-2.0*\Offset,\Height,-1*\Offset + \Depth);
  \draw[white, thick, fill=myblue] (lbm0) -- (lbm1) -- (lbm2) -- (lbm3) -- cycle;
  \draw[white, thick] ($(lbm0)!0.5!(lbm3)$) -- ($(lbm1)!0.5!(lbm2)$);

  \coordinate (lbo0) at (-3.5*\Offset,0,-1*\Offset + \Depth);
  \coordinate (lbo1) at (-4.5*\Offset,0,-1*\Offset + \Depth);
  \coordinate (lbo2) at (-4.5*\Offset,\Height,-1*\Offset + \Depth);
  \coordinate (lbo3) at (-3.5*\Offset,\Height,-1*\Offset + \Depth);
  \draw[white, thick, fill=myblue] (lbo0) -- (lbo1) -- (lbo2) -- (lbo3) -- cycle;
  \draw[white, thick] ($(lbo0)!0.5!(lbo3)$) -- ($(lbo1)!0.5!(lbo2)$);

  \draw[white, thick, fill=myred] (lbi0) -- (lbi1) --
  ($(lbi1)!0.5!(lbi2)$) -- ($(lbi0)!0.5!(lbi3)$) -- (lbi0) -- cycle;
  \draw[white, thick, fill=myred] (lbm0) -- (lbm1) --
  ($(lbm1)!0.5!(lbm2)$) -- ($(lbm0)!0.5!(lbm3)$) -- (lbm0) -- cycle;
  \draw[white, thick, fill=myred] (lbo0) -- (lbo1) --
  ($(lbo1)!0.5!(lbo2)$) -- ($(lbo0)!0.5!(lbo3)$) -- (lbo0) -- cycle;

  \coordinate (lfi0) at (-0.5*\Offset,0,\Depth);
  \coordinate (lfi1) at (-1.5*\Offset,0,\Depth);
  \coordinate (lfi2) at (-1.5*\Offset,\Height,\Depth);
  \coordinate (lfi3) at (-0.5*\Offset,\Height,\Depth);
  \draw[white, thick, fill=myblue] (lfi0) -- (lfi1) -- (lfi2) -- (lfi3) -- cycle;
  \draw[white, thick] ($(lfi0)!0.5!(lfi3)$) -- ($(lfi1)!0.5!(lfi2)$);

  \coordinate (lfm0) at (-2.0*\Offset,0,\Depth);
  \coordinate (lfm1) at (-3.0*\Offset,0,\Depth);
  \coordinate (lfm2) at (-3.0*\Offset,\Height,\Depth);
  \coordinate (lfm3) at (-2.0*\Offset,\Height,\Depth);
  \draw[white, thick, fill=myblue] (lfm0) -- (lfm1) -- (lfm2) -- (lfm3) -- cycle;
  \draw[white, thick] ($(lfm0)!0.5!(lfm3)$) -- ($(lfm1)!0.5!(lfm2)$);

  \coordinate (lfo0) at (-3.5*\Offset,0,\Depth);
  \coordinate (lfo1) at (-4.5*\Offset,0,\Depth);
  \coordinate (lfo2) at (-4.5*\Offset,\Height,\Depth);
  \coordinate (lfo3) at (-3.5*\Offset,\Height,\Depth);
  \draw[white, thick, fill=myblue] (lfo0) -- (lfo1) -- (lfo2) -- (lfo3) -- cycle;
  \draw[white, thick] ($(lfo0)!0.5!(lfo3)$) -- ($(lfo1)!0.5!(lfo2)$);

  \draw[white, thick, fill=myred] (lfi0) -- (lfi1) --
  ($(lfi1)!0.5!(lfi2)$) -- ($(lfi0)!0.5!(lfi3)$) -- (lfi0) -- cycle;
  \draw[white, thick, fill=myred] (lfm0) -- (lfm1) --
  ($(lfm1)!0.5!(lfm2)$) -- ($(lfm0)!0.5!(lfm3)$) -- (lfm0) -- cycle;
  \draw[white, thick, fill=myred] (lfo0) -- (lfo1) --
  ($(lfo1)!0.5!(lfo2)$) -- ($(lfo0)!0.5!(lfo3)$) -- (lfo0) -- cycle;

  \coordinate (bfu0) at (0,-1.5*\Offset,\Depth);
  \coordinate (bfu1) at (0,-0.5*\Offset,\Depth);
  \coordinate (bfu2) at (\Width,-0.5*\Offset,\Depth);
  \coordinate (bfu3) at (\Width,-1.5*\Offset,\Depth);

  \coordinate (bbu0) at (0,-1.5*\Offset,\Depth - \Offset);
  \coordinate (bbu1) at (0,-0.5*\Offset,\Depth - \Offset);
  \coordinate (bbu2) at (\Width,-0.5*\Offset,\Depth - \Offset);
  \coordinate (bbu3) at (\Width,-1.5*\Offset,\Depth - \Offset);

  \coordinate (bfl0) at (0,-3.0*\Offset,\Depth);
  \coordinate (bfl1) at (0,-2.0*\Offset,\Depth);
  \coordinate (bfl2) at (\Width,-2.0*\Offset,\Depth);
  \coordinate (bfl3) at (\Width,-3.0*\Offset,\Depth);

  \coordinate (bbl0) at (0,-3.0*\Offset,\Depth - \Offset);
  \coordinate (bbl1) at (0,-2.0*\Offset,\Depth -1*\Offset);
  \coordinate (bbl2) at (\Width,-2.0*\Offset,\Depth - \Offset);
  \coordinate (bbl3) at (\Width,-3.0*\Offset,\Depth - \Offset);

  \draw[white, thick, fill=myblue] (bbu1) -- (bbu2) -- (bbu3) -- (bbu0) -- cycle;

  \draw[white, thick] ($(bbu0)!.33!(bbu3)$) -- ($(bbu1)!.33!(bbu2)$);
  \draw[white, thick] ($(bbu0)!.67!(bbu3)$) -- ($(bbu1)!.67!(bbu2)$);

  \draw[white, thick, fill=myblue] (bfu1) -- (bfu2) -- (bfu3) -- (bfu0) -- cycle;

  \draw[white, thick] ($(bfu0)!.33!(bfu3)$) -- ($(bfu1)!.33!(bfu2)$);
  \draw[white, thick] ($(bfu0)!.67!(bfu3)$) -- ($(bfu1)!.67!(bfu2)$);

  \draw[white, thick, fill=myblue] (bbl1) -- (bbl2) -- (bbl3) -- (bbl0) -- cycle;
  \draw[white, thick] ($(bbl0)!.33!(bbl3)$) -- ($(bbl1)!.33!(bbl2)$);
  \draw[white, thick] ($(bbl0)!.67!(bbl3)$) -- ($(bbl1)!.67!(bbl2)$);

  \draw[white, thick, fill=myblue] (bfl1) -- (bfl2) -- (bfl3) -- (bfl0) -- cycle;
  \draw[white, thick] ($(bfl0)!.33!(bfl3)$) -- ($(bfl1)!.33!(bfl2)$);
  \draw[white, thick] ($(bfl0)!.67!(bfl3)$) -- ($(bfl1)!.67!(bfl2)$);

  \draw[white, thick, fill=myorange] (bbu3) -- ($(bbu0)!.33!(bbu3)$) --
  ($(bbu1)!.33!(bbu2)$) -- (bbu2) -- (bbu3) -- cycle;
  \draw[white, thick, fill=myorange] (bbl3) -- ($(bbl0)!.33!(bbl3)$) --
  ($(bbl1)!.33!(bbl2)$) -- (bbl2) -- (bbl3) -- cycle;
  \draw[white, thick, fill=myorange] (bfu3) -- ($(bfu0)!.33!(bfu3)$) --
  ($(bfu1)!.33!(bfu2)$) -- (bfu2) -- (bfu3) -- cycle;
  \draw[white, thick, fill=myorange] (bfl3) -- ($(bfl0)!.33!(bfl3)$) --
  ($(bfl1)!.33!(bfl2)$) -- (bfl2) -- (bfl3) -- cycle;

  \draw[white, thick, fill=myred] (bbu3) -- ($(bbu0)!.67!(bbu3)$) --
  ($(bbu1)!.67!(bbu2)$) -- (bbu2) -- (bbu3) -- cycle;
  \draw[white, thick, fill=myred] (bbl3) -- ($(bbl0)!.67!(bbl3)$) --
  ($(bbl1)!.67!(bbl2)$) -- (bbl2) -- (bbl3) -- cycle;
  \draw[white, thick, fill=myred] (bfu3) -- ($(bfu0)!.67!(bfu3)$) --
  ($(bfu1)!.67!(bfu2)$) -- (bfu2) -- (bfu3) -- cycle;
  \draw[white, thick, fill=myred] (bfl3) -- ($(bfl0)!.67!(bfl3)$) --
  ($(bfl1)!.67!(bfl2)$) -- (bfl2) -- (bfl3) -- cycle;

  \draw[<->, thick] ($(lfo1) - (0.5*\Offset, 0, 0)$) -- ($(lfo2) - (0.5*\Offset,
  0, 0)$)  node[midway, rotate=90, above] {Distributed};

  \draw[<->, thick] ($(lbo2) + (0.1*\Offset,0.5*\Offset,0)$) -- ($(lbi3)+ (0,0.5*\Offset,0)$);
  \draw[<->, thick] ($(lfo2)+ (0,0.5*\Offset,0)$) -- ($(lbo2) + (0.0*\Offset,0.5*\Offset,0)$) node [above] {Replicated};

  \node (xlabel) at (1.0\Width, 1.5\Height, \Depth) {$\tns{X}$};
  \node (alabel) at ($(lfo0)!0.5!(lfi1) - (0,0.5*\Offset,0)$) {$\mat{A}$};
  \node (blabel) at ($(bfl0)!0.5!(bfl3) - (0,0.5*\Offset,0)$) {$\mat{B}$};
  \node (clabel) at ($(rli0)!0.5!(rlo3) - (0,0.5*\Offset,0)$) {$\mat{C}$};

\end{tikzpicture}

%% file: images/genten_two_sided.tex
\newcommand{\Depth}{2}
\newcommand{\Height}{2}
\newcommand{\Width}{2}
\newcommand{\Offset}{0.5}
\begin{tikzpicture}

  \coordinate (origin) at (0,0,0);
  \coordinate (lowerbackright) at (\Width,0,0);
  \coordinate (upperbackright) at (\Width,\Height,0);
  \coordinate (upperbackleft) at (0,\Height,0);

  \coordinate (lowerfrontleft) at (0,0,\Depth);
  \coordinate (lowerfrontright) at (\Width,0,\Depth);
  \coordinate (upperfrontright) at (\Width,\Height,\Depth);
  \coordinate (upperfrontleft) at (0,\Height,\Depth);

  \coordinate(rui0) at (\Width+0.5*\Offset,0+0.5*\Height,\Depth);
  \coordinate(rui1) at (\Width+0.5*\Offset,0+0.5*\Height,0);
  \coordinate(rui2) at (\Width+1.5*\Offset,0+0.5*\Height,0);
  \coordinate(rui3) at (\Width+1.5*\Offset,0+0.5*\Height,\Depth);
  \draw[white, thick, fill=white] (rui0) -- (rui1) -- (rui2) -- (rui3) -- cycle;
  \draw[white, thick] ($(rui0)!0.5!(rui1)$) -- ($(rui2)!0.5!(rui3)$);

  \coordinate(rli0) at (\Width+0.5*\Offset,0,\Depth);
  \coordinate(rli1) at (\Width+0.5*\Offset,0,0);
  \coordinate(rli2) at (\Width+1.5*\Offset,0,0);
  \coordinate(rli3) at (\Width+1.5*\Offset,0,\Depth);
  \draw[white, thick, fill=white] (rli0) -- (rli1) -- (rli2) -- (rli3) -- cycle;
  \draw[white, thick] ($(rli0)!0.5!(rli1)$) -- ($(rli2)!0.5!(rli3)$);

  \coordinate(rum0) at (\Width+2.0*\Offset,0+0.5*\Height,\Depth);
  \coordinate(rum1) at (\Width+2.0*\Offset,0+0.5*\Height,0);
  \coordinate(rum2) at (\Width+3.0*\Offset,0+0.5*\Height,0);
  \coordinate(rum3) at (\Width+3.0*\Offset,0+0.5*\Height,\Depth);
  \draw[white, thick, fill=white] (rum0) -- (rum1) -- (rum2) -- (rum3) -- cycle;
  \draw[white, thick] ($(rum0)!0.5!(rum1)$) -- ($(rum2)!0.5!(rum3)$);

  \coordinate(rlm0) at (\Width+2.0*\Offset,0,\Depth);
  \coordinate(rlm1) at (\Width+2.0*\Offset,0,0);
  \coordinate(rlm2) at (\Width+3.0*\Offset,0,0);
  \coordinate(rlm3) at (\Width+3.0*\Offset,0,\Depth);
  \draw[white, thick, fill=white] (rlm0) -- (rlm1) -- (rlm2) -- (rlm3) -- cycle;
  \draw[white, thick] ($(rlm0)!0.5!(rlm1)$) -- ($(rlm2)!0.5!(rlm3)$);

  \coordinate(ruo0) at (\Width+3.5*\Offset,0+0.5*\Height,\Depth);
  \coordinate(ruo1) at (\Width+3.5*\Offset,0+0.5*\Height,0);
  \coordinate(ruo2) at (\Width+4.5*\Offset,0+0.5*\Height,0);
  \coordinate(ruo3) at (\Width+4.5*\Offset,0+0.5*\Height,\Depth);
  \draw[white, thick, fill=white] (ruo0) -- (ruo1) -- (ruo2) -- (ruo3) -- cycle;
  \draw[white, thick] ($(ruo0)!0.5!(ruo1)$) -- ($(ruo2)!0.5!(ruo3)$);

  \coordinate(rlo0) at (\Width+3.5*\Offset,0,\Depth);
  \coordinate(rlo1) at (\Width+3.5*\Offset,0,0);
  \coordinate(rlo2) at (\Width+4.5*\Offset,0,0);
  \coordinate(rlo3) at (\Width+4.5*\Offset,0,\Depth);
  \draw[white, thick, fill=white] (rlo0) -- (rlo1) -- (rlo2) -- (rlo3) -- cycle;
  \draw[white, thick] ($(rlo0)!0.5!(rlo1)$) -- ($(rlo2)!0.5!(rlo3)$);

  \draw[white, thick, fill=white] (rlo0) -- ($(rlo0)!0.5!(rlo1)$) --
  ($(rlo2)!0.5!(rlo3)$) -- (rlo3) -- (rlo0) -- cycle;

  \coordinate (lbi0) at (-0.5*\Offset,0,-1*\Offset + \Depth);
  \coordinate (lbi1) at (-1.5*\Offset,0,-1*\Offset + \Depth);
  \coordinate (lbi2) at (-1.5*\Offset,\Height,-1*\Offset + \Depth);
  \coordinate (lbi3) at (-0.5*\Offset,\Height,-1*\Offset + \Depth);
  \draw[white, thick, fill=white] (lbi0) -- (lbi1) -- (lbi2) -- (lbi3) -- cycle;
  \draw[white, thick] ($(lbi0)!0.5!(lbi3)$) -- ($(lbi1)!0.5!(lbi2)$);

  \coordinate (lbm0) at (-2.0*\Offset,0,-1*\Offset + \Depth);
  \coordinate (lbm1) at (-3.0*\Offset,0,-1*\Offset + \Depth);
  \coordinate (lbm2) at (-3.0*\Offset,\Height,-1*\Offset + \Depth);
  \coordinate (lbm3) at (-2.0*\Offset,\Height,-1*\Offset + \Depth);
  \draw[white, thick, fill=white] (lbm0) -- (lbm1) -- (lbm2) -- (lbm3) -- cycle;
  \draw[white, thick] ($(lbm0)!0.5!(lbm3)$) -- ($(lbm1)!0.5!(lbm2)$);

  \coordinate (lbo0) at (-3.5*\Offset,0,-1*\Offset + \Depth);
  \coordinate (lbo1) at (-4.5*\Offset,0,-1*\Offset + \Depth);
  \coordinate (lbo2) at (-4.5*\Offset,\Height,-1*\Offset + \Depth);
  \coordinate (lbo3) at (-3.5*\Offset,\Height,-1*\Offset + \Depth);
  \draw[white, thick, fill=white] (lbo0) -- (lbo1) -- (lbo2) -- (lbo3) -- cycle;
  \draw[white, thick] ($(lbo0)!0.5!(lbo3)$) -- ($(lbo1)!0.5!(lbo2)$);

  \coordinate (lfi0) at (-0.5*\Offset,0,\Depth);
  \coordinate (lfi1) at (-1.5*\Offset,0,\Depth);
  \coordinate (lfi2) at (-1.5*\Offset,\Height,\Depth);
  \coordinate (lfi3) at (-0.5*\Offset,\Height,\Depth);
  \draw[white, thick, fill=white] (lfi0) -- (lfi1) -- (lfi2) -- (lfi3) -- cycle;
  \draw[white, thick] ($(lfi0)!0.5!(lfi3)$) -- ($(lfi1)!0.5!(lfi2)$);

  \coordinate (lfm0) at (-2.0*\Offset,0,\Depth);
  \coordinate (lfm1) at (-3.0*\Offset,0,\Depth);
  \coordinate (lfm2) at (-3.0*\Offset,\Height,\Depth);
  \coordinate (lfm3) at (-2.0*\Offset,\Height,\Depth);
  \draw[white, thick, fill=white] (lfm0) -- (lfm1) -- (lfm2) -- (lfm3) -- cycle;
  \draw[white, thick] ($(lfm0)!0.5!(lfm3)$) -- ($(lfm1)!0.5!(lfm2)$);

  \coordinate (lfo0) at (-3.5*\Offset,0,\Depth);
  \coordinate (lfo1) at (-4.5*\Offset,0,\Depth);
  \coordinate (lfo2) at (-4.5*\Offset,\Height,\Depth);
  \coordinate (lfo3) at (-3.5*\Offset,\Height,\Depth);
  \draw[white, thick, fill=white] (lfo0) -- (lfo1) -- (lfo2) -- (lfo3) -- cycle;
  \draw[white, thick] ($(lfo0)!0.5!(lfo3)$) -- ($(lfo1)!0.5!(lfo2)$);

  \draw[white, thick, fill=teal] (lfi0) -- (lfi1) --
  ($(lfi1)!0.5!(lfi2)$) -- ($(lfi0)!0.5!(lfi3)$) -- (lfi0) -- cycle;

  \coordinate (bfu0) at (0,-1.5*\Offset,\Depth);
  \coordinate (bfu1) at (0,-0.5*\Offset,\Depth);
  \coordinate (bfu2) at (\Width,-0.5*\Offset,\Depth);
  \coordinate (bfu3) at (\Width,-1.5*\Offset,\Depth);

  \coordinate (bbu0) at (0,-1.5*\Offset,\Depth - \Offset);
  \coordinate (bbu1) at (0,-0.5*\Offset,\Depth - \Offset);
  \coordinate (bbu2) at (\Width,-0.5*\Offset,\Depth - \Offset);
  \coordinate (bbu3) at (\Width,-1.5*\Offset,\Depth - \Offset);

  \coordinate (bfl0) at (0,-3.0*\Offset,\Depth);
  \coordinate (bfl1) at (0,-2.0*\Offset,\Depth);
  \coordinate (bfl2) at (\Width,-2.0*\Offset,\Depth);
  \coordinate (bfl3) at (\Width,-3.0*\Offset,\Depth);

  \coordinate (bbl0) at (0,-3.0*\Offset,\Depth - \Offset);
  \coordinate (bbl1) at (0,-2.0*\Offset,\Depth -1*\Offset);
  \coordinate (bbl2) at (\Width,-2.0*\Offset,\Depth - \Offset);
  \coordinate (bbl3) at (\Width,-3.0*\Offset,\Depth - \Offset);

  \draw[white, thick, fill=white] (bbu1) -- (bbu2) -- (bbu3) -- (bbu0) -- cycle;

  \draw[white, thick] ($(bbu0)!.33!(bbu3)$) -- ($(bbu1)!.33!(bbu2)$);
  \draw[white, thick] ($(bbu0)!.67!(bbu3)$) -- ($(bbu1)!.67!(bbu2)$);

  \draw[white, thick, fill=white] (bfu1) -- (bfu2) -- (bfu3) -- (bfu0) -- cycle;

  \draw[white, thick] ($(bfu0)!.33!(bfu3)$) -- ($(bfu1)!.33!(bfu2)$);
  \draw[white, thick] ($(bfu0)!.67!(bfu3)$) -- ($(bfu1)!.67!(bfu2)$);

  \draw[white, thick, fill=white] (bbl1) -- (bbl2) -- (bbl3) -- (bbl0) -- cycle;
  \draw[white, thick] ($(bbl0)!.33!(bbl3)$) -- ($(bbl1)!.33!(bbl2)$);
  \draw[white, thick] ($(bbl0)!.67!(bbl3)$) -- ($(bbl1)!.67!(bbl2)$);

  \draw[white, thick, fill=white] (bfl1) -- (bfl2) -- (bfl3) -- (bfl0) -- cycle;
  \draw[white, thick] ($(bfl0)!.33!(bfl3)$) -- ($(bfl1)!.33!(bfl2)$);
  \draw[white, thick] ($(bfl0)!.67!(bfl3)$) -- ($(bfl1)!.67!(bfl2)$);

  \draw[white, thick, fill=white] (bfu3) -- ($(bfu0)!.67!(bfu3)$) --
  ($(bfu1)!.67!(bfu2)$) -- (bfu2) -- (bfu3) -- cycle;
  \draw[white, thick, fill=white] (bfl3) -- ($(bfl0)!.67!(bfl3)$) --
  ($(bfl1)!.67!(bfl2)$) -- (bfl2) -- (bfl3) -- cycle;

  \draw[<->, thick, white] ($(lfo1) - (0.5*\Offset, 0, 0)$) -- ($(lfo2) - (0.5*\Offset,
  0, 0)$)  node[midway, rotate=90, above] {Distributed};

  \draw[<->, thick, white] ($(lbo2) + (0.1*\Offset,0.5*\Offset,0)$) -- ($(lbi3)+ (0,0.5*\Offset,0)$);
  \draw[<->, thick, white] ($(lfo2)+ (0,0.5*\Offset,0)$) -- ($(lbo2) + (0.0*\Offset,0.5*\Offset,0)$) node [above] {Replicated};

  \node (alabel) at ($(lfo0)!0.5!(lfi1) - (0,0.5*\Offset,0)$) [white] {$\mat{A}$};
  \node (blabel) at ($(bfl0)!0.5!(bfl3) - (0,0.5*\Offset,0)$) [white]  {$\mat{B}$};
  \node (clabel) at ($(rli0)!0.5!(rlo3) - (0,0.5*\Offset,0)$) [white] {$\mat{C}$};

  \input{images/tensor3d.tex}

  \coordinate(rli0) at (\Width+0.5*\Offset,0,\Depth);
  \coordinate(rli1) at (\Width+0.5*\Offset,0,0);
  \coordinate(rli2) at (\Width+1.5*\Offset,0,0);
  \coordinate(rli3) at (\Width+1.5*\Offset,0,\Depth);

  \draw[white, semithick, fill=myblue] (rli0) -- ($(rli0)!1.000!(rli1)$) -- ($(rli3)!1.000!(rli2)$) -- (rli3) -- cycle;
  \draw[white, semithick, fill=myblue] (rli0) -- ($(rli0)!0.916!(rli1)$) -- ($(rli3)!0.916!(rli2)$) -- (rli3) -- cycle;
  \draw[white, semithick, fill=myblue] (rli0) -- ($(rli0)!0.833!(rli1)$) -- ($(rli3)!0.833!(rli2)$) -- (rli3) -- cycle;
  \draw[white, semithick, fill=myblue] (rli0) -- ($(rli0)!0.750!(rli1)$) -- ($(rli3)!0.750!(rli2)$) -- (rli3) -- cycle;
  \draw[white, semithick, fill=myblue] (rli0) -- ($(rli0)!0.666!(rli1)$) -- ($(rli3)!0.666!(rli2)$) -- (rli3) -- cycle;
  \draw[white, semithick, fill=myblue] (rli0) -- ($(rli0)!0.583!(rli1)$) -- ($(rli3)!0.583!(rli2)$) -- (rli3) -- cycle;
  \draw[white, semithick, fill=myred] (rli0) -- ($(rli0)!0.500!(rli1)$) -- ($(rli3)!0.500!(rli2)$) -- (rli3) -- cycle;
  \draw[white, semithick, fill=myred] (rli0) -- ($(rli0)!0.416!(rli1)$) -- ($(rli3)!0.416!(rli2)$) -- (rli3) -- cycle;
  \draw[white, semithick, fill=myred] (rli0) -- ($(rli0)!0.333!(rli1)$) -- ($(rli3)!0.333!(rli2)$) -- (rli3) -- cycle;
  \draw[white, semithick, fill=myred] (rli0) -- ($(rli0)!0.250!(rli1)$) -- ($(rli3)!0.250!(rli2)$) -- (rli3) -- cycle;
  \draw[white, semithick, fill=myred] (rli0) -- ($(rli0)!0.167!(rli1)$) -- ($(rli3)!0.167!(rli2)$) -- (rli3) -- cycle;
  \draw[white, semithick, fill=myred] (rli0) -- ($(rli0)!0.083!(rli1)$) -- ($(rli3)!0.083!(rli2)$) -- (rli3) -- cycle;

  \coordinate (lfi0) at (-0.5*\Offset,0,\Depth);
  \coordinate (lfi1) at (-1.5*\Offset,0,\Depth);
  \coordinate (lfi2) at (-1.5*\Offset,\Height,\Depth);
  \coordinate (lfi3) at (-0.5*\Offset,\Height,\Depth);
  \draw[white, thick, fill=myblue] (lfi0) -- (lfi1) -- (lfi2) -- (lfi3) -- cycle;

  \draw[white, thick, fill=myblue] (lfi0) -- (lfi1) -- ($(lfi1)!1.000!(lfi2)$) -- ($(lfi0)!1.000!(lfi3)$) -- (lfi0) -- cycle;
  \draw[white, thick, fill=myblue] (lfi0) -- (lfi1) -- ($(lfi1)!0.916!(lfi2)$) -- ($(lfi0)!0.916!(lfi3)$) -- (lfi0) -- cycle;
  \draw[white, thick, fill=myblue] (lfi0) -- (lfi1) -- ($(lfi1)!0.833!(lfi2)$) -- ($(lfi0)!0.833!(lfi3)$) -- (lfi0) -- cycle;
  \draw[white, thick, fill=myblue] (lfi0) -- (lfi1) -- ($(lfi1)!0.750!(lfi2)$) -- ($(lfi0)!0.750!(lfi3)$) -- (lfi0) -- cycle;
  \draw[white, thick, fill=myblue] (lfi0) -- (lfi1) -- ($(lfi1)!0.666!(lfi2)$) -- ($(lfi0)!0.666!(lfi3)$) -- (lfi0) -- cycle;
  \draw[white, thick, fill=myblue] (lfi0) -- (lfi1) -- ($(lfi1)!0.583!(lfi2)$) -- ($(lfi0)!0.583!(lfi3)$) -- (lfi0) -- cycle;
  \draw[white, thick, fill=myred] (lfi0) -- (lfi1) -- ($(lfi1)!0.500!(lfi2)$) -- ($(lfi0)!0.500!(lfi3)$) -- (lfi0) -- cycle;
  \draw[white, thick, fill=myred] (lfi0) -- (lfi1) -- ($(lfi1)!0.416!(lfi2)$) -- ($(lfi0)!0.416!(lfi3)$) -- (lfi0) -- cycle;
  \draw[white, thick, fill=myred] (lfi0) -- (lfi1) -- ($(lfi1)!0.333!(lfi2)$) -- ($(lfi0)!0.333!(lfi3)$) -- (lfi0) -- cycle;
  \draw[white, thick, fill=myred] (lfi0) -- (lfi1) -- ($(lfi1)!0.250!(lfi2)$) -- ($(lfi0)!0.250!(lfi3)$) -- (lfi0) -- cycle;
  \draw[white, thick, fill=myred] (lfi0) -- (lfi1) -- ($(lfi1)!0.167!(lfi2)$) -- ($(lfi0)!0.167!(lfi3)$) -- (lfi0) -- cycle;
  \draw[white, thick, fill=myred] (lfi0) -- (lfi1) -- ($(lfi1)!0.083!(lfi2)$) -- ($(lfi0)!0.083!(lfi3)$) -- (lfi0) -- cycle;

  \coordinate (bfu0) at (0,-1.5*\Offset,\Depth);
  \coordinate (bfu1) at (0,-0.5*\Offset,\Depth);
  \coordinate (bfu2) at (\Width,-0.5*\Offset,\Depth);
  \coordinate (bfu3) at (\Width,-1.5*\Offset,\Depth);

  \draw[white, thick, fill=myblue] (bfu1) -- (bfu2) -- (bfu3) -- (bfu0) -- cycle;

  \draw[white, thick,    fill=myred] (bfu0) -- (bfu1) -- ($(bfu1)!1.000!(bfu2)$) -- ($(bfu0)!1.000!(bfu3)$) -- (bfu0) -- cycle;
  \draw[white, thick,    fill=myred] (bfu0) -- (bfu1) -- ($(bfu1)!0.916!(bfu2)$) -- ($(bfu0)!0.916!(bfu3)$) -- (bfu0) -- cycle;
  \draw[white, thick,    fill=myred] (bfu0) -- (bfu1) -- ($(bfu1)!0.833!(bfu2)$) -- ($(bfu0)!0.833!(bfu3)$) -- (bfu0) -- cycle;
  \draw[white, thick,    fill=myred] (bfu0) -- (bfu1) -- ($(bfu1)!0.750!(bfu2)$) -- ($(bfu0)!0.750!(bfu3)$) -- (bfu0) -- cycle;
  \draw[white, thick, fill=myorange] (bfu0) -- (bfu1) -- ($(bfu1)!0.666!(bfu2)$) -- ($(bfu0)!0.666!(bfu3)$) -- (bfu0) -- cycle;
  \draw[white, thick, fill=myorange] (bfu0) -- (bfu1) -- ($(bfu1)!0.583!(bfu2)$) -- ($(bfu0)!0.583!(bfu3)$) -- (bfu0) -- cycle;
  \draw[white, thick, fill=myorange] (bfu0) -- (bfu1) -- ($(bfu1)!0.500!(bfu2)$) -- ($(bfu0)!0.500!(bfu3)$) -- (bfu0) -- cycle;
  \draw[white, thick, fill=myorange] (bfu0) -- (bfu1) -- ($(bfu1)!0.416!(bfu2)$) -- ($(bfu0)!0.416!(bfu3)$) -- (bfu0) -- cycle;
  \draw[white, thick,  fill=myblue] (bfu0) -- (bfu1) -- ($(bfu1)!0.333!(bfu2)$) -- ($(bfu0)!0.333!(bfu3)$) -- (bfu0) -- cycle;
  \draw[white, thick,  fill=myblue] (bfu0) -- (bfu1) -- ($(bfu1)!0.250!(bfu2)$) -- ($(bfu0)!0.250!(bfu3)$) -- (bfu0) -- cycle;
  \draw[white, thick,  fill=myblue] (bfu0) -- (bfu1) -- ($(bfu1)!0.167!(bfu2)$) -- ($(bfu0)!0.167!(bfu3)$) -- (bfu0) -- cycle;
  \draw[white, thick,  fill=myblue] (bfu0) -- (bfu1) -- ($(bfu1)!0.083!(bfu2)$) -- ($(bfu0)!0.083!(bfu3)$) -- (bfu0) -- cycle;

  \draw[] ($(lfi1) - (0.5*\Offset, 0, 0)$) -- ($(lfi2) - (0.5*\Offset,
  0, 0)$)  node[midway, rotate=90, above] {Distributed};

  \draw[<->, thick] ($(lfi1) - (0.5*\Offset, 0, 0)$) -- ($(lfi2) - (0.5*\Offset,
  0.5*\Height, 0)$)  node[midway, rotate=90, above] {};
  \draw[<->, thick] ($(lfi2) - (0.5*\Offset, 0.5*\Height, 0)$) -- ($(lfi2) - (0.5*\Offset,
  0, 0)$)  node[midway, rotate=90, above] {};

  \node (xlabel) at (1.0\Width, 1.5\Height, \Depth) {$\tns{X}$};
  \node (alabel) at ($(lfi0)!0.5!(lfi1) - (0,0.5*\Offset,0)$) {$\mat{A}$};
  \node (blabel) at ($(bfu0)!0.5!(bfu3) - (0,0.5*\Offset,0)$) {$\mat{B}$};
  \node (clabel) at ($(rli0)!0.5!(rli3) - (0,0.5*\Offset,0)$) {$\mat{C}$};

\end{tikzpicture}

%% file: fed-opt.tex
{\gcp}-{\fedadam} is an extension of Reddi et al.'s
{\fedopt}~\cite{Reddi21AdaptiveFederatedOptimization} to GCP.
{\gcp}-{\fedadam} is presented in \Cref{alg:gcp-fed-iter}. The
{\fedopt} step
(Lines~\ref{alg:gcp-fed-iter:diff}--\ref{alg:gcp-fed-iter:assign})
incorporates the progress made locally
(Lines~\ref{alg:gcp-fed-iter:gradient}--\ref{alg:gcp-fed-iter:update_local})
into the global model. The tensor $\tns{D}$ is a finite difference
approximation of the true gradient. The call to \texttt{.update()}
(Line~\ref{alg:gcp-fed-iter:update_global};
i.e.,~\cref{alg:adam-update}) computes the bias-corrected moment
estimates from the current model and gradient estimate to update the
model parameters, i.e., the {\adam} step in~\cite[Alg.
1]{Kingma15AdamMethodStochastic}.

\begin{algorithm}
  \caption{\textsc{FedAdam} mini-batch iteration (asynchronous)}
  \label{alg:gcp-fed-iter}
  \begin{algorithmic}[1]
    \Function{FedAdamIteration}{Global model $\tns{M}$, local model
    $\tns{U}$, $sampler$, epoch $e > 0$, asynchrony $\tau > 0$}
    \If{$\text{mod}(e, \tau) = 0$}
    \State $\tns{D} \gets \tns{U} - \tns{M}$ \label{alg:gcp-fed-iter:diff}
    \State $AllReduce(\tns{D})$ \label{alg:gcp-fed-iter:allreduce}
    \State $\tns{U}.update(\tns{D})$ \label{alg:gcp-fed-iter:update_global}
    \State $\tns{M} \gets \tns{U}$ \label{alg:gcp-fed-iter:assign}
    \EndIf
    \State $\tns{G} \gets sampler.Gradient(\tns{M})$
    \label{alg:gcp-fed-iter:gradient}
    \State $\tns{M}.update(\tns{G})$ \label{alg:gcp-fed-iter:update_local}
    \EndFunction
  \end{algorithmic}
\end{algorithm}%

%% file: experiments.tex

\subsection{Synchronous Parallel GCP Performance}

\begin{table}[ht]
  \centering
  \caption{Data tensors used in the numerical experiments.  The arXiv
    data set was first studied in
    \cite{Phipps23StreamingGeneralizedCanonicala}.  The Amazon data set
  is from FROSTT~\cite{frosttdataset}.}
  \label{tab:tensors}
  \begin{tabular}{lcrc}
    \toprule
    Tensor       & Dimensions  & Nonzeros  & Loss model  \\
    \midrule
    Synthetic & $300 \times 200 \times 100$ & $59\text{K}$ & Poisson \\
    arXiv  & $172 \times 25\text{K} \times 300$  & 30M  & Poisson \\
    Amazon & $4.8\text{M} \times 1.7\text{M} \times 1.8\text{M}$ &
    ~1.7B & Poisson \\
    \bottomrule
  \end{tabular}
\end{table}

\begin{figure}
  \centering
  \begin{subfigure}[t]{0.45\textwidth}
    \centering
    \begin{adjustbox}{width=\linewidth}
      \includegraphics[width=\textwidth]{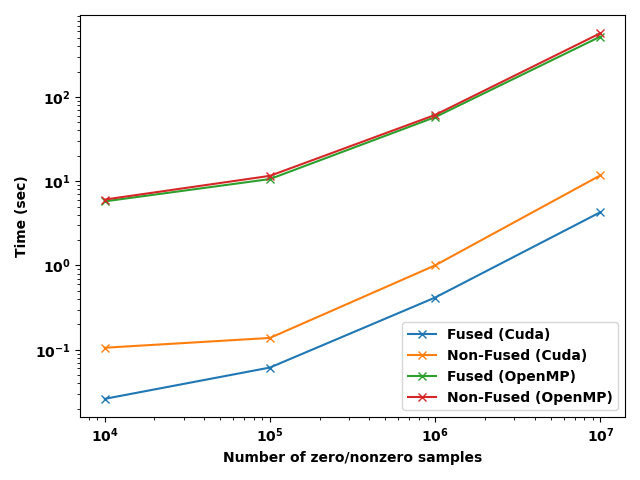}
    \end{adjustbox}
    \caption{Run time}
    \label{fig:fused-perf-time}
  \end{subfigure}
  \begin{subfigure}[t]{0.45\textwidth}
    \centering
    \begin{adjustbox}{width=\linewidth}
      \includegraphics[width=\textwidth]{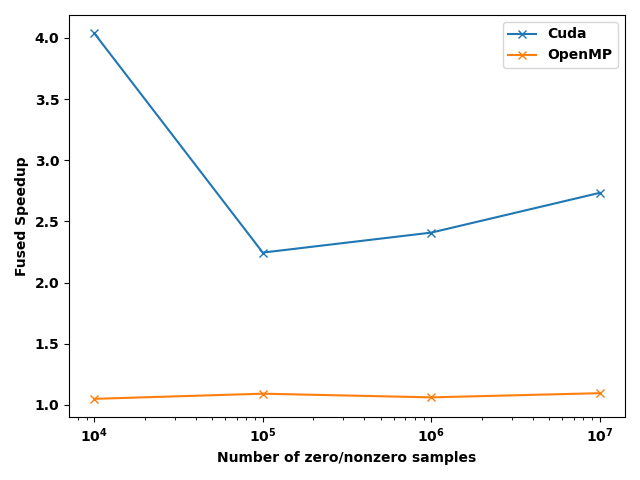}
    \end{adjustbox}
    \caption{Fused speedup}
    \label{fig:fused-perf-speedup}
  \end{subfigure}
  \caption{Comparison of fused and non-fused semi-stratified sampling/MTTKRP
    kernels on a single node with CUDA and OpenMP backends for the arXiv tensor.
    Stochastic gradient evaluation time
    (\subref{fig:fused-perf-time}) and speedup
    of fused over non-fused (\subref{fig:fused-perf-speedup}), indicating 2-4x
  speedup on the GPU and about 10\% speedup the CPU.}
  \label{fig:fused-perf}
\end{figure}

In this section we investigate the performance of distributed and shared memory
parallelism for GCP based on the synchronous algorithm described by
\cref{alg:gcp-sgd-iter}.  Since shared-memory parallelism with Kokkos for MTTKRP
was extensively discussed in~\cite{Phipps19SoftwareSparseTensora}, we only
consider portions specific to GCP, in particular comparing fused and non-fused
sampling/MTTKRP, and the two distributed parallelism approaches described above.
The computer hardware for this evaluation is a single 18-node rack within
a larger computing cluster where each node consists of one 32-core AMD EPYC
7543P processor, 512 GB of DDR4 memory, 1 NVIDIA A100 GPU (40 GB of HBM2) and
a 100 Gb/s 18-way flat Infiniband interconnect.     We compiled GenTen with GCC
13.3.0, NVCC 12.6.2, and OpenMPI 4.1.7, enabling three choices of Kokkos
back-ends:  Serial, OpenMP, and CUDA.  We consider GCP decomposition of the
tensors listed in \cref{tab:tensors} in three configurations:  MPI+CUDA with
1 MPI process per node using CUDA for the A100 GPU, MPI+OpenMP with 1 MPI
process per node and 32 OpenMP threads per processes (each bound to a core), and
MPI+Serial with 32 MPI processes per node (each bound to a core).  First,
\cref{fig:fused-perf} compares the runtime cost of the sampled gradient
evaluation for the arXiv tensor on a single compute node over two epochs (with
100 iterations per epoch) using the fused and non-fused sampling/MTTKRP
approaches on both the CPU and GPU architectures with varying numbers of
semi-stratified gradient samples per iteration (with the same number of zero and
nonzero samples).  On the GPU, we use GenTen's atomic-based MTTKRP/fused
implementations and because of the small mode sizes, the thread-privatization
with OpenMP on the CPU.   We use 1.5M zero and nonzero samples for the objective
function evaluation and default values for the {\adam} step.  We observe a 2-4x
reduction in cost for the fused approach on the GPU architecture and about
a 10\% reduction on the CPU.

\begin{figure}
  \centering
  \begin{subfigure}[t]{0.45\textwidth}
    \centering
    \begin{adjustbox}{width=\linewidth}
      \includegraphics[width=\textwidth]{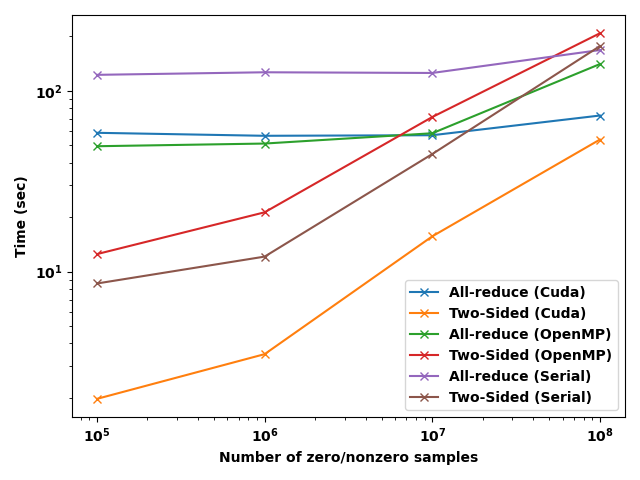}
    \end{adjustbox}
    \caption{\textit{Distributed GCP runtime cost comparing
    all-reduce and two-sided with Serial, OpenMP, and CUDA backends.}}
    \label{fig:distributed-gcp-cost}
  \end{subfigure}\hfill
  \begin{subfigure}[t]{0.45\textwidth}
    \centering
    \begin{adjustbox}{width=\linewidth}
      \includegraphics[width=\textwidth]{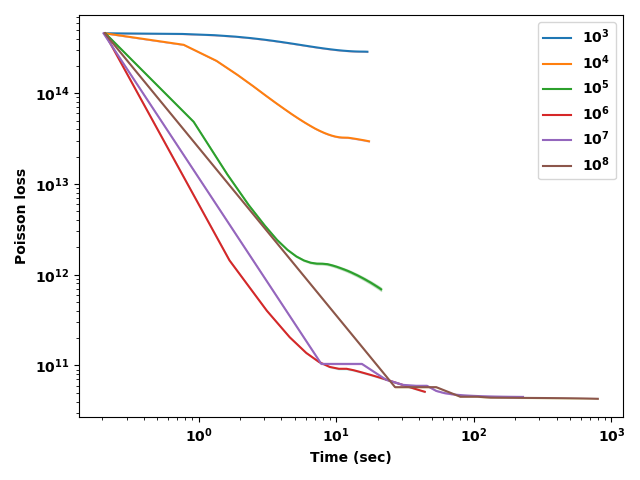}
    \end{adjustbox}
    \caption{\textit{GCP-Amazon convergence with varying numbers of
    zero/nonzero samples.}}
    \label{fig:distributed-gcp-convergence}
  \end{subfigure}
  \caption{Comparison of distributed GCP approaches for over a range
  of total number of samples.}
  \label{fig:distributed-gcp-perf}
\end{figure}

Next, we consider the all-reduce and two-sided distributed parallelism
approaches described in \cref{sec:gcp_mpi} using the Amazon tensor distributed
across 18 compute nodes.  \Cref{fig:distributed-gcp-cost} displays the total GCP
runtime cost of two GCP epochs (with 100 iterations per epoch) for both
approaches on the CPU and GPU architectures with varying numbers of
semi-stratified zero/nonzero samples per iteration (with 10M zero/nonzero
objective function samples and default {\adam} parameters).  As the two-sided
approach requires the sampled tensor indices in order to construct the parallel
communication patterns between processors, it cannot use the fused approach.  We
therefore use non-fused for both approaches to isolate the parallel
communication costs.  In this case however, we use atomic-based MTTKRP on both
the GPU with CUDA and the CPU with OpenMP (but GenTen's serial implementation
for MPI-only) due to the large tensor mode sizes.  We see, as expected, the cost
of the all-reduce approach is essentially flat as the number of samples
increases, since it performs a dense communication after each MTTKRP calculation
independent of the number of samples, whereas the cost of the two-sided approach
grows substantially.  However, different cross-over points are observed
depending on the chosen architecture, i.e., about $10^8$ samples for MPI+CUDA
and MPI-only, and about $10^7$ samples for MPI+OpenMP.  For MPI+CUDA, this is
likely due to the increased cost of the sends and receives directly from/to the
GPU as required by the two-sided approach.  However, different numbers of
samples also result in different time-to-solution since convergence is typically
quicker (in terms of iterations) with more samples.  This is investigated in
\cref{fig:distributed-gcp-convergence} by displaying the computed Poisson loss
averaged over 5 trials as a function of computing time for a range of
zero/nonzero samples using the two-sided approach and MPI+CUDA, indicating
a range of $10^6$--$10^8$ samples is optimal.  In this range, two-sided is
generally more efficient, but not always.  Hence we see that the best choice of
distributed MTTKRP method depends on both the architecture and the number of
chosen samples (which is problem dependent), with either the all-reduce or
two-sided approaches being more efficient.

\begin{figure}
  \centering
  \begin{subfigure}[t]{0.45\textwidth}
    \centering
    \begin{adjustbox}{width=\linewidth}
      \includegraphics[width=\textwidth]{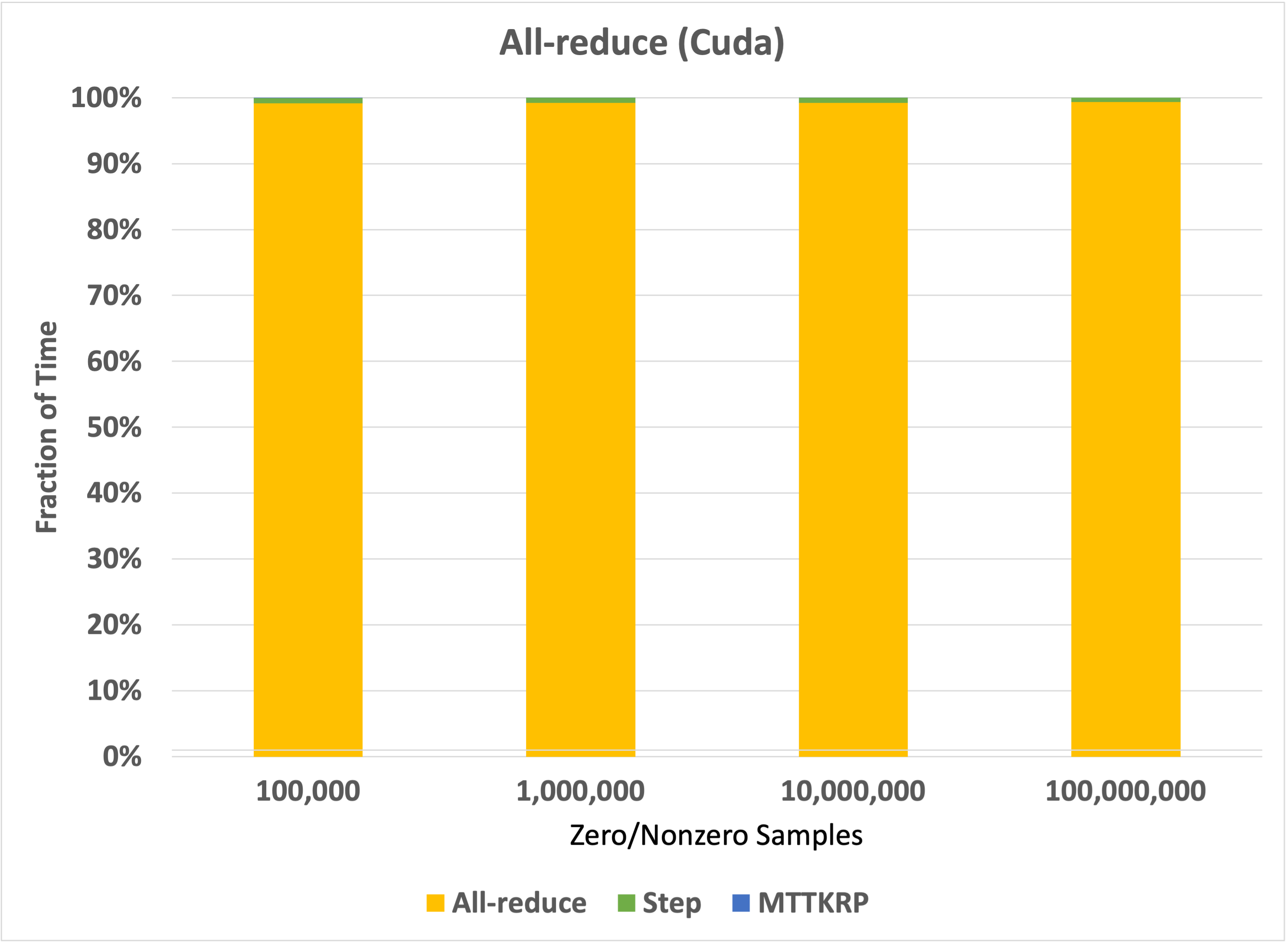}
    \end{adjustbox}
    \caption{}
  \end{subfigure}\hfill
  \begin{subfigure}[t]{0.45\textwidth}
    \centering
    \begin{adjustbox}{width=\linewidth}
      \includegraphics[width=\textwidth]{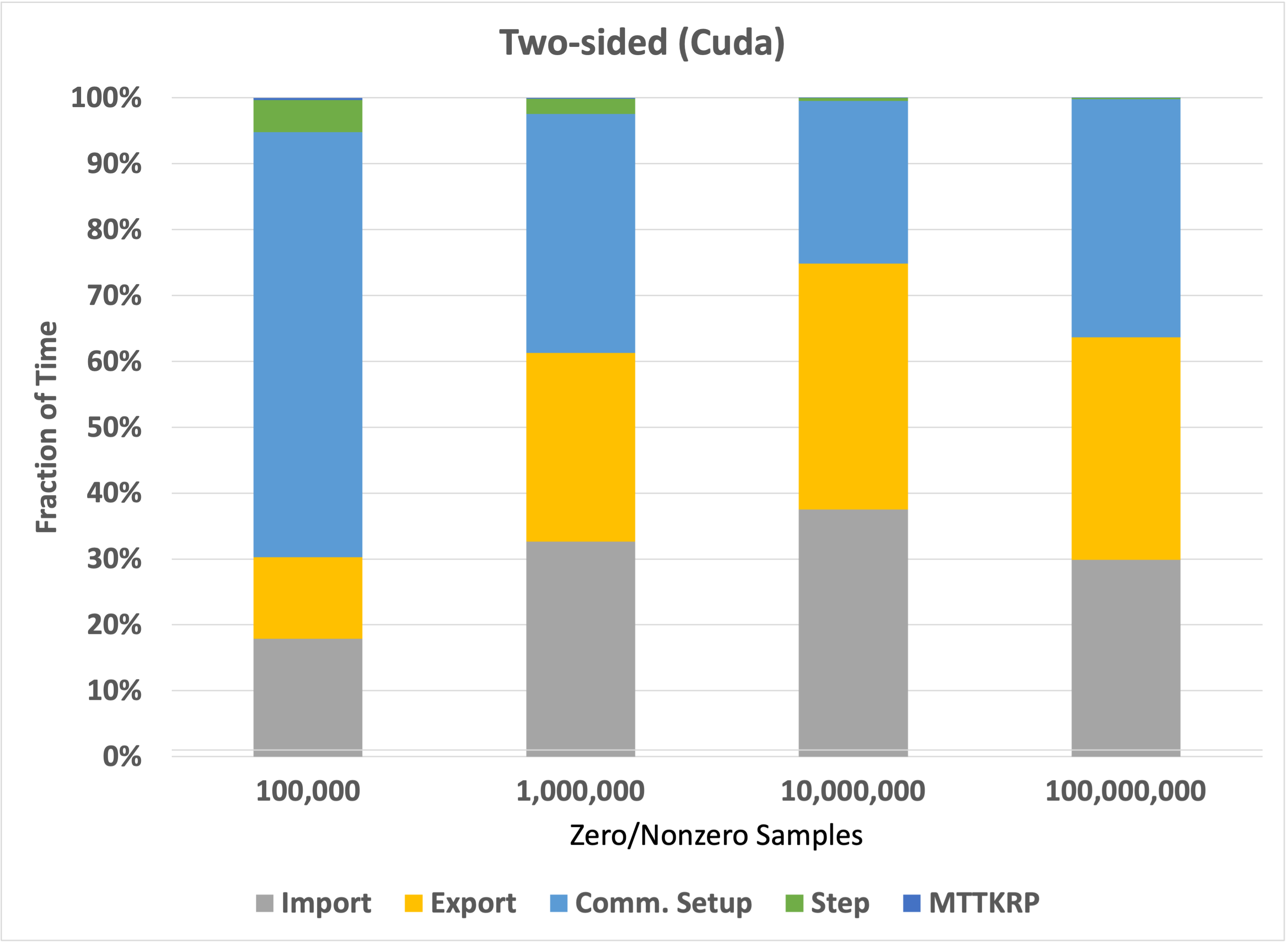}
    \end{adjustbox}
    \caption{}
  \end{subfigure}
  \begin{subfigure}[t]{0.45\textwidth}
    \centering
    \begin{adjustbox}{width=\linewidth}
      \includegraphics[width=\textwidth]{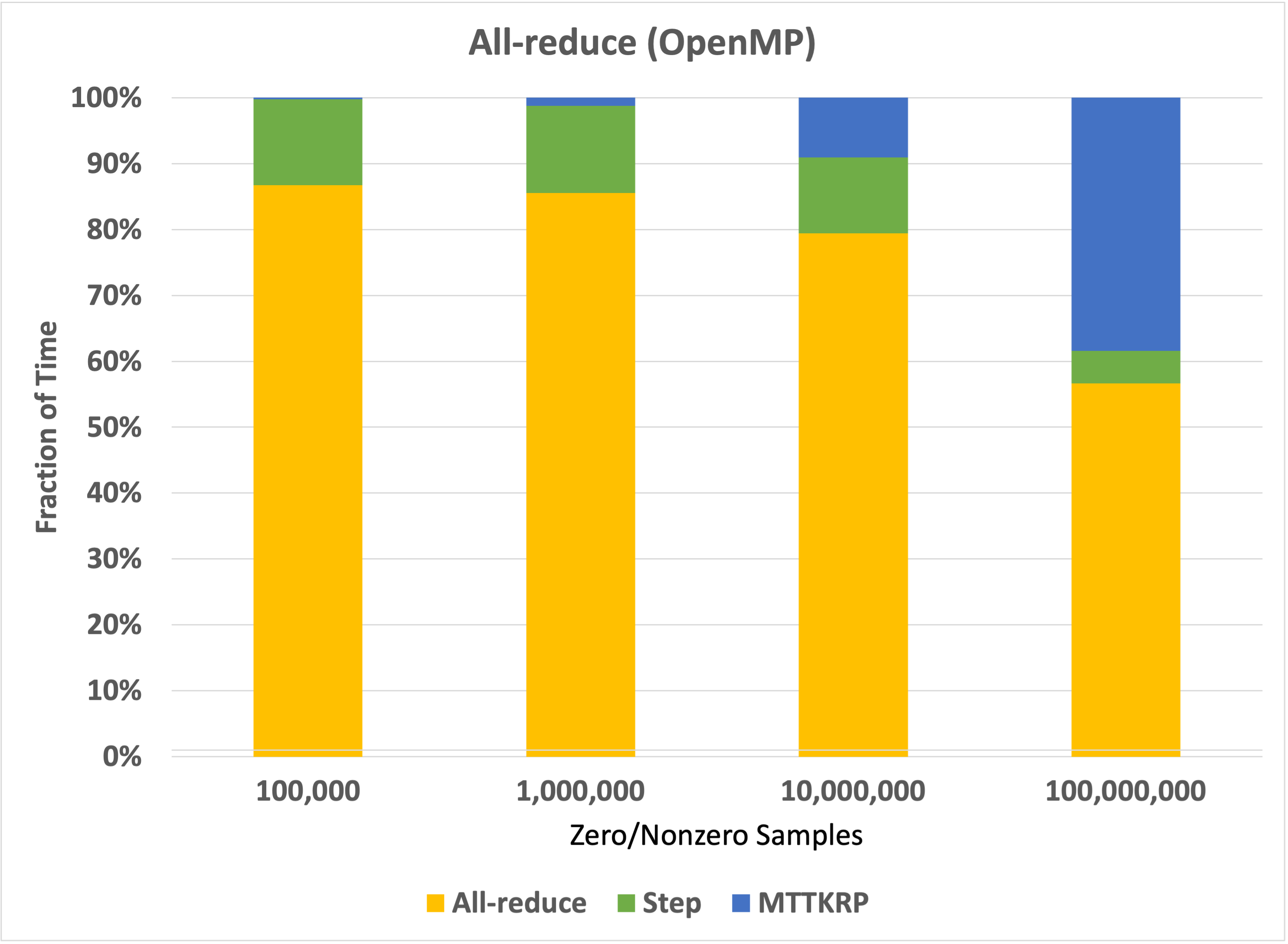}
    \end{adjustbox}
    \caption{}
  \end{subfigure}\hfill
  \begin{subfigure}[t]{0.45\textwidth}
    \centering
    \begin{adjustbox}{width=\linewidth}
      \includegraphics[width=\textwidth]{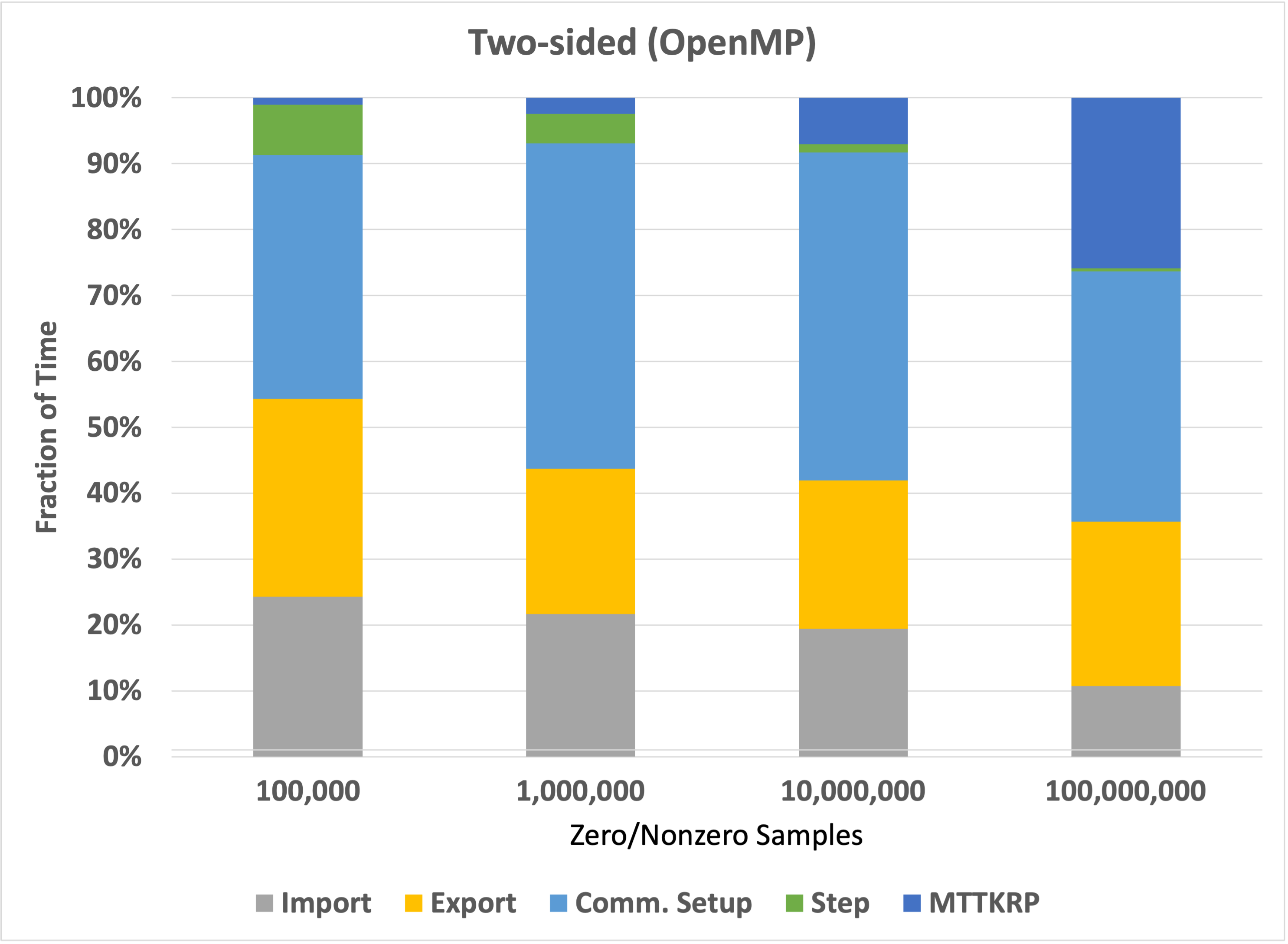}
    \end{adjustbox}
    \caption{}
  \end{subfigure}
  \begin{subfigure}[t]{0.45\textwidth}
    \centering
    \begin{adjustbox}{width=\linewidth}
      \includegraphics[width=\textwidth]{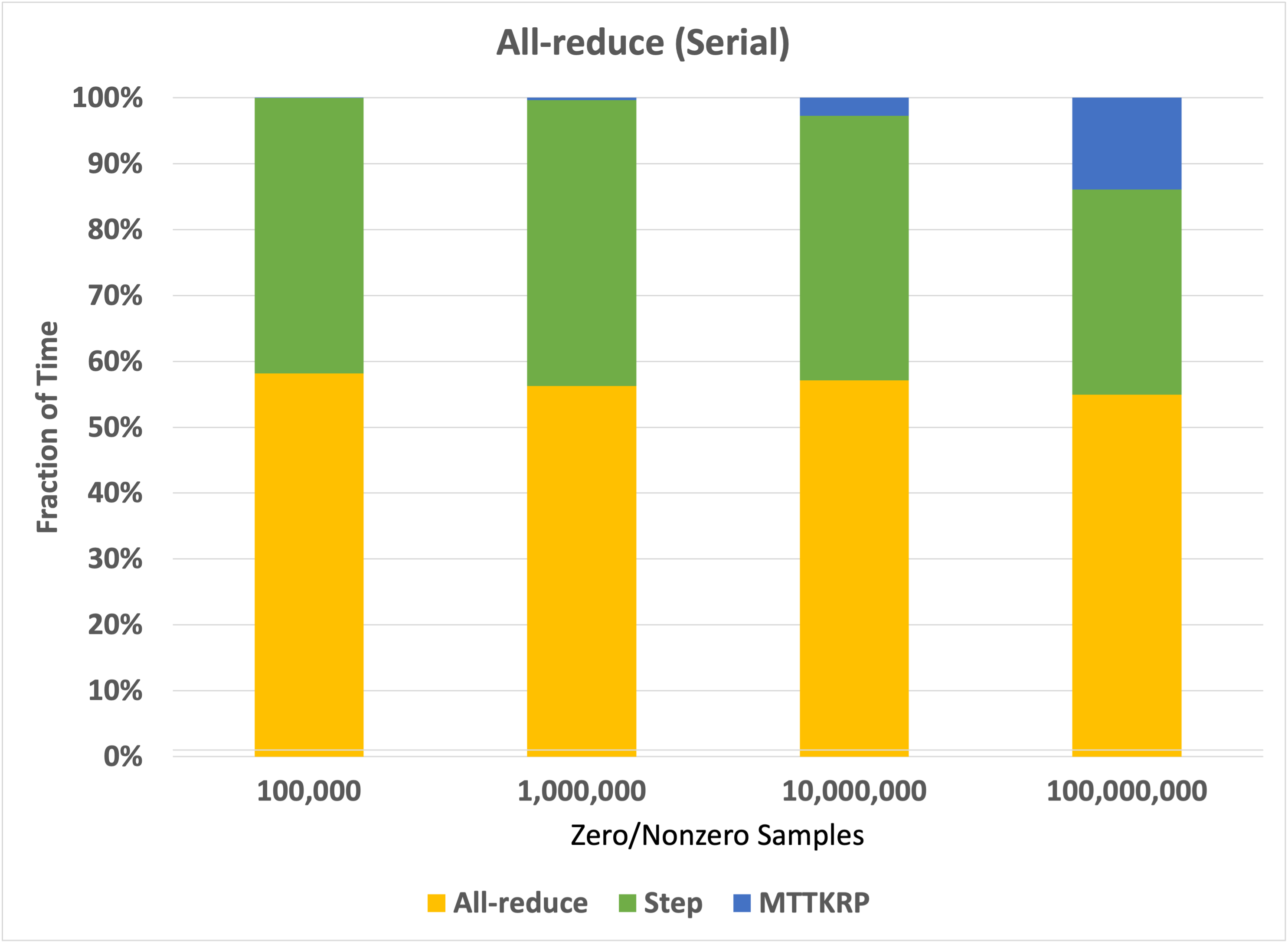}
    \end{adjustbox}
    \caption{}
  \end{subfigure}\hfill
  \begin{subfigure}[t]{0.45\textwidth}
    \centering
    \begin{adjustbox}{width=\linewidth}
      \includegraphics[width=\textwidth]{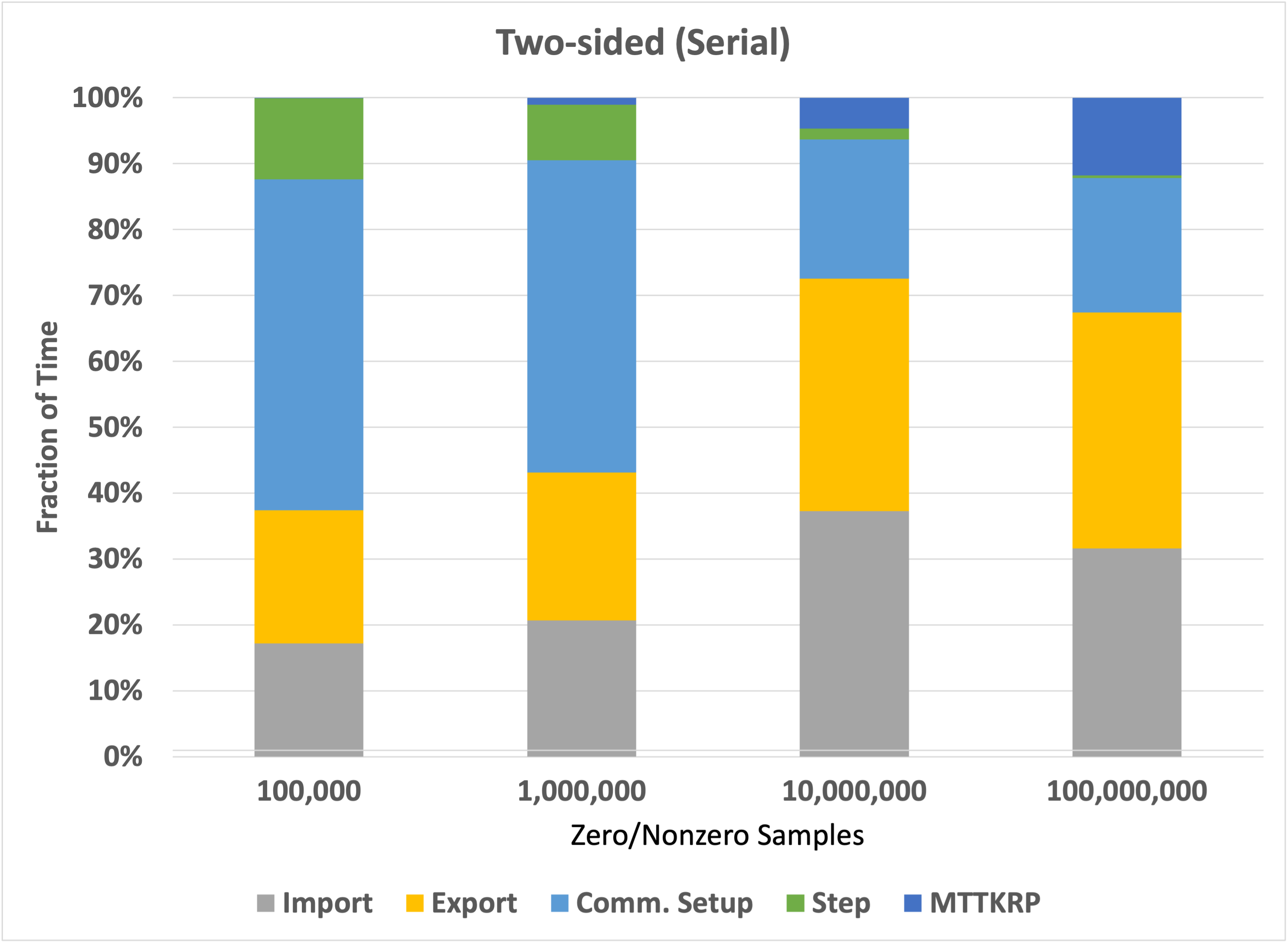}
    \end{adjustbox}
    \caption{}
  \end{subfigure}
  \caption{Breaking down GCP runtime for the all-reduce and two-sided
  communication approaches on CUDA, OpenMP, and Serial.}
  \label{fig:distributed-gcp-perf-fractions}
\end{figure}

To analyze the performance of the two parallel distributed approaches, fractions
of the total {\gcp} runtime for each number of samples and Kokkos backend are
shown in \cref{fig:distributed-gcp-perf-fractions} in terms of parallel
communication (All-reduce for the all-reduce method or factor matrix import,
export, and communication setup for the two-sided method), the {\adam} step, and
local MTTKRP.  We see that for the all-reduce method, the parallel all-reduce
required is a substantial fraction of the runtime in all cases.  However the
local MTTKRP cost is also significant in the OpenMP case with larger numbers of
samples, likely due to increased frequency of thread collisions in the atomic
MTTKRP update.  We also see the {\adam} step cost is significant for the
MPI-only Serial case due to the reduced parallelism from the replicated factor
matrices.  For the two-sided approach, communication costs dominate in all
cases, in particular the setup costs from determining the communication pattern
after each sampled tensor are substantial.  Thus a method that could eliminate
these setup costs would be beneficial.

\begin{figure}
  \centering
  \begin{adjustbox}{width=0.5\linewidth}
    \includegraphics[width=\textwidth]{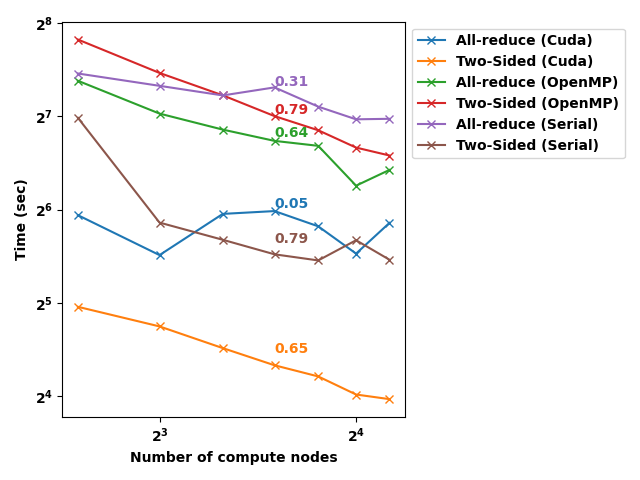}
  \end{adjustbox}
  \caption{Strong scaling of two GCP-{\adam} epochs across a range of
    compute nodes.  The slope of the best linear fit (in log-space) is
  shown above each curve.}
  \label{fig:gcp-strong}
\end{figure}

Finally, strong scaling of the above experiment from 6-18 compute nodes (since
  6 A100 GPUs is the fewest that will fit the Amazon tensor and its GCP
decomposition) using $10^7$ zero/nonzero samples is shown in
\cref{fig:gcp-strong}.  We see nearly linear scaling for two-sided and
all-reduce with OpenMP, but much poorer scaling for all-reduce, in particular
essentially no speedup on the GPU with all-reduce.

\subsection{Convergence implications for synchronous and asynchronous
approaches}
In this section, we describe numerical experiments comparing
{\gcpsgd} and {\gcpfed} to understand the implication of asynchrony
on GCP convergence and performance.
First, we performed a parameter sensitivity analysis on synthetic data utilizing
on-node parallelism only to establish baseline behavior.
Second, we introduced distributed parallelism and conducted experiments with two
large-scale, sparse tensor datasets.

Since federated
approaches may promote rapid progress on ``easier'' subsets of heterogeneous
data, diminishing the contributions learned on more important subsets, we felt
that the homogeneity across clients resulting from synthetic data generation
made a true distributed (i.e., multiple node) study unwarranted. Instead, we
emphasize the impact of the finite-difference approximation in {\gcpfed} on
convergence, which is absent from {\gcpsgd}.

\subsubsection{Hyperparameter study on synthetic data}
A key challenge of such a
comparison is that the set of algorithm hyperparameters that balances
these objectives may differ between algorithms. Hyperparameter tuning
is an oft-derided but important task for both algorithm convergence
and computational performance.

Rather than exhaustive grid search, which grows exponentially in the
number of hyperparameters and their values, we pursue Latin hypercube
sampling (LHS)~\cite{McKay79ComparisonThreeMethods} for
hyperparameter tuning. LHS generates sampling distributions that are
useful in assessing parameter uncertainties in a ``global'' sense.
LHS partitions the parameter space in such a way that samples are
chosen to guarantee every row and column in the hypercube has exactly
one sample. We use the Dakota~\cite{DakotaMultilevelParallel20} to
generate LHS samples,
schedule jobs (i.e., trials), and we leverage the provided
functionality that allows us to easily identify impactful
hyperparameters and their values based on a rank-ordering of their
Spearman rank correlation coefficients (also known as simple rank
correlation coefficients).

\paragraph{Experiment design}
Our hyperparameter study proceeds as follows in a two-stage process.
In the first stage, $s$ samples are drawn uniformly for
each parameter between upper and lower bounds. An LHS sample is thus
a tuple with each element a value corresponding to a specific
parameter. For each LHS sample parameterization, we generate $n$ Kruskal tensor
starting points and run {\gcpsgd} and {\gcpfed} from each initial
guess. Both methods also share the random seed used by GenTen to sample the
tensor for loss value estimation and pseudo-gradient computations. The final
loss values are then recorded for both methods. Once the first stage
completes, we search for the minimum across all $2sn$ instances and
compute the relative errors with respect to the minimum. We
return the errors to Dakota, which reports how each parameterization
correlates to convergence to the maximum likelihood estimator. Specifically, we
use Spearman's rank correlation (also referred to as \textit{simple} rank
correlation) to measure the strength of monotonic relationships between
hyperparameter values and their impact on the final loss. The utility of using
simple rank correlations is that it allows one to measure sensitivity to the
order of magnitude of a hyperparameter value, rather than a specific value. In
this way, we hope that our conclusions will retain generality for both our
experiments specifically and for stochastic GCP solvers more broadly.

We pursue a refinement policy for the two-stage parameter study.
After the second stage completes, we set $s \gets 2s$, restart the
process, and continue until a sufficient number of samples have been
collected. In this way, we can observe convergence behavior in the
correlation coefficients.

In our hyperparameter study, we used $n=10$ starting points and
initialized the refinement policy with $s=14$ samples.\footnote{The Dakota
  User's Guide recommends $2d$ samples at a minimum, where $d$ is the number of
hyperparameters under study.} The parameters of interest, their upper and lower
bounds, and values or ranges of values from the prior literature are reported
in~\Cref{tab:lhs:poisson:continuous}. Note that the \textit{meta-rate} and
\textit{downpour-iterations} parameters are used for with {\gcpfed} only.

We initially conducted our study on two synthetic datasets with samples randomly
chosen from Gaussian and Poisson distributions and ran both methods with the
respective loss function. However, we observed that the convergence
results for the synthetic
Gaussian tensor depended entirely on the initial guess and random seed---i.e.,
both methods converged to the maximum likelihood estimator with the same
empirical probability irrespective of the hyperparameterization---so
we omitted results for this dataset.

\begin{table}[!ht]
  \centering
  \caption{Continuous uniform uncertain parameters for initial LHS study.}
  \label{tab:lhs:poisson:continuous}
  \begin{tabular}{lrrr}
    \toprule
    Parameter & Lower Bound & Upper Bound & Default \\
    \midrule
    \textit{rate} & $10^{-5.5}$ & $10^{-0.5}$ & $10^{-3}$
    \cite{Kingma15AdamMethodStochastic}\\
    \textit{decay} & $10^{-3.5}$ & $10^{-0.5}$ & $10^{-1}$
    \cite{Kingma15AdamMethodStochastic}, $[10^{-3},10^{-1}]$
    \cite{Reddi21AdaptiveFederatedOptimization}\\
    \textit{adam-beta1} & 0.0 & 1.0 & 0.9 \cite{Kingma15AdamMethodStochastic}\\
    \textit{adam-beta2} & 0.9 & 0.99999 & 0.999
    \cite{Kingma15AdamMethodStochastic}\\
    \textit{adam-eps} & $10^{-16}$ & $10^{-10}$ & $10^{-8}$
    \cite{Kingma15AdamMethodStochastic},  $10^{-3}$
    \cite{Reddi21AdaptiveFederatedOptimization}\\
    \textit{meta-rate}\textsuperscript{\textdagger} & $10^{-5.5}$ &
    $10^{-0.5}$ & $[10^{-3}, 10^{-0.5}]$
    \cite{Reddi21AdaptiveFederatedOptimization}\\
    \textit{downpour-iterations}\textsuperscript{\textdagger} & 1 &
    256 & 1 \cite{Dean12LargeScaleDistributed}, 10
    \cite{Zhang15DeepLearningElastic} \\
    \bottomrule
  \end{tabular}\\
  \textsuperscript{\textdagger} {\gcpfed} only
\end{table}

\begin{figure}
  \centering
  \begin{subfigure}[t]{0.5\textwidth}
    \centering
    \includegraphics[width=\textwidth]{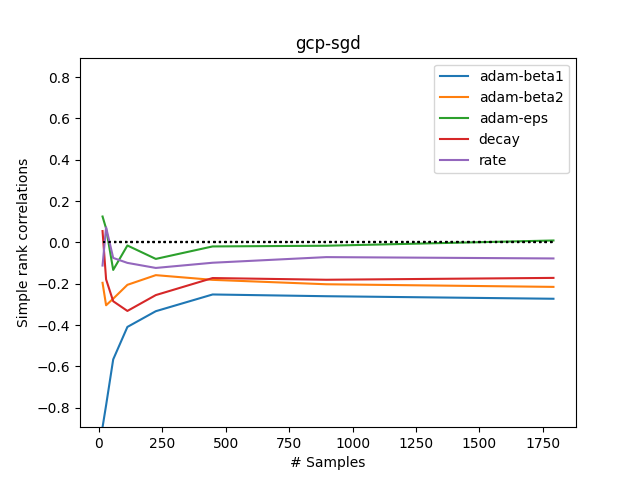}
    \caption{{\gcpsgd}}
    \label{fig:experiments:correlations:gcp-sgd}
  \end{subfigure}\hfill
  \begin{subfigure}[t]{0.5\textwidth}
    \centering
    \includegraphics[width=\textwidth]{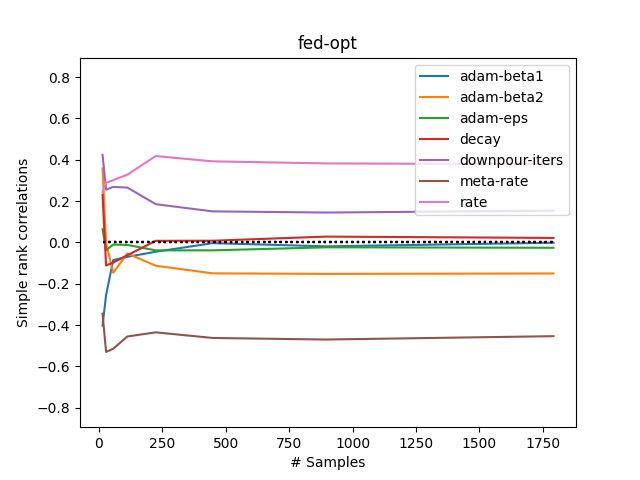}
    \caption{{\gcpfed}}
    \label{fig:experiments:correlations:fed-opt}
  \end{subfigure}
  \caption{Simple rank correlation coefficient convergence of {\gcpsgd} (left)
    and {\gcpfed} (right) for the synthetic count dataset as the number of LHS
  samples grows.}
  \label{fig:experiments:correlations}
\end{figure}

\paragraph{Results}
We display the convergence of the simple rank correlation coefficients from the
parameter study for the synthetic count dataset
in~\cref{fig:experiments:correlations}. A useful heuristic is that stronger
correlations are represented by more extreme values along the $y$-axis. For
example, \textit{adam-beta1} is more strongly correlated to final loss value
than \textit{adam-eps} for {\gcpsgd}. These figures help visualize a rank
ordering by parameter impact and are useful for identifying which parameters
warrant further study. In fact, these results are somewhat surprising for
{\gcpsgd}, since we anticipated the learning rate parameters \textit{rate} and
\textit{decay} to dominate the simple rank correlation coefficients. By
contrast, the results for {\gcpfed}, where the server and client learning rate
parameters are most impactful, are closer to our expectations.

\begin{figure}
  \centering
  \includegraphics[width=\textwidth]{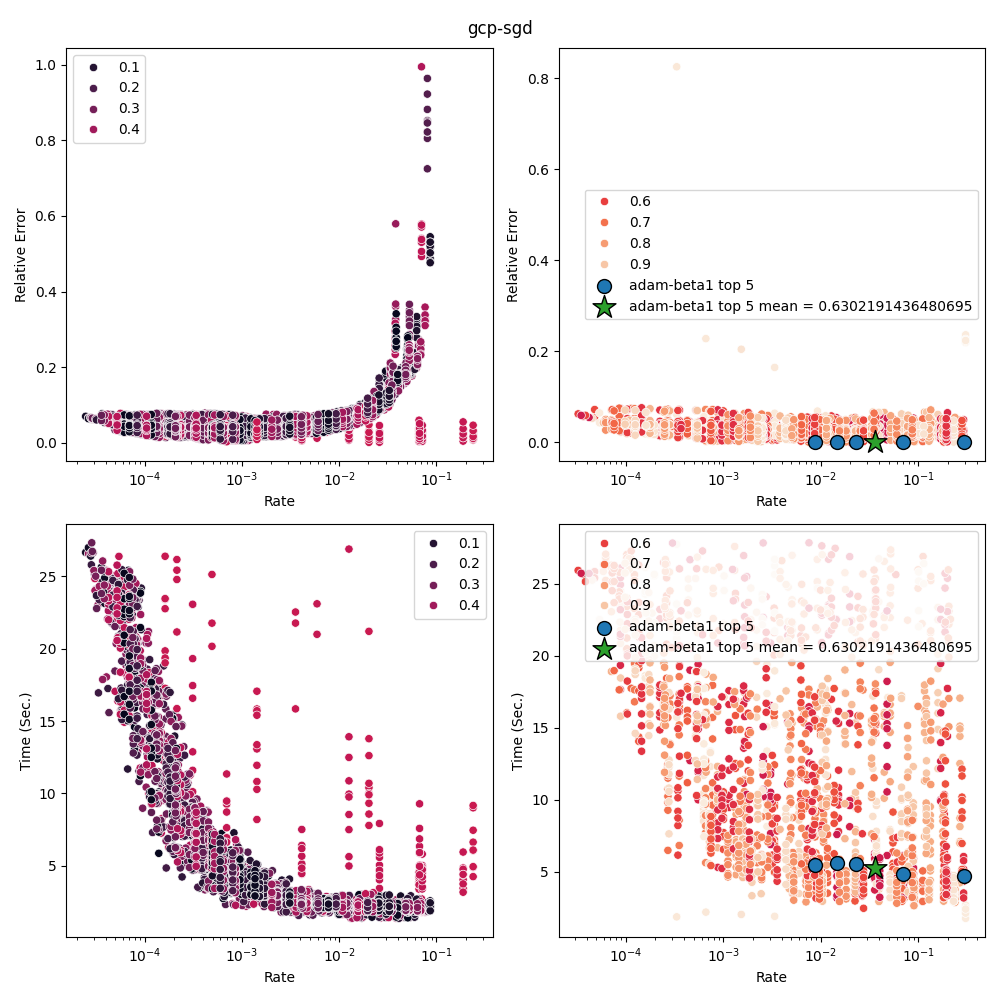}
  \caption{Scatter plots varying \textit{rate} and \textit{adam-beta1} for
    {\gcpsgd}. The $y$-axis denotes the relative error from the MLE (top) and
    total time (bottom); \textit{rate} is plotted along the $x$-axis;
    variations in
  \textit{adam-beta1} are denoted by the hue.}
  \label{fig:experiments:dakota-poisson:gcp-sgd:rate-adam-beta1}
\end{figure}

\begin{figure}
  \centering
  \includegraphics[width=\textwidth]{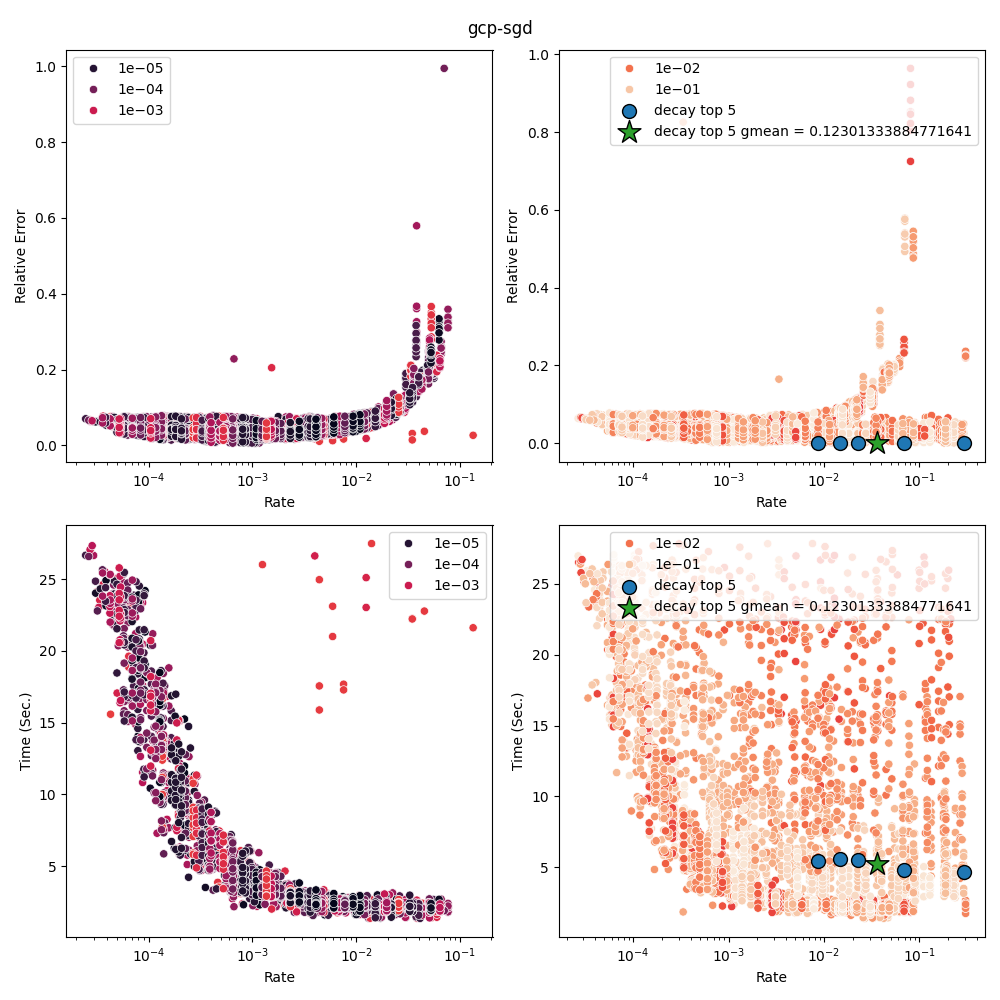}
  \caption{Scatter plots varying \textit{rate} and \textit{decay} for {\gcpsgd}.
    The $y$-axis denotes the relative error from the MLE (top) and total time
    (bottom); \textit{rate} is plotted along the $x$-axis; variations in
  \textit{decay} are denoted by the hue.}
  \label{fig:experiments:dakota-poisson:gcp-sgd:rate-decay}
\end{figure}

\begin{figure}
  \centering
  \begin{subfigure}[t]{0.85\textwidth}
    \centering
    \includegraphics[width=\textwidth]{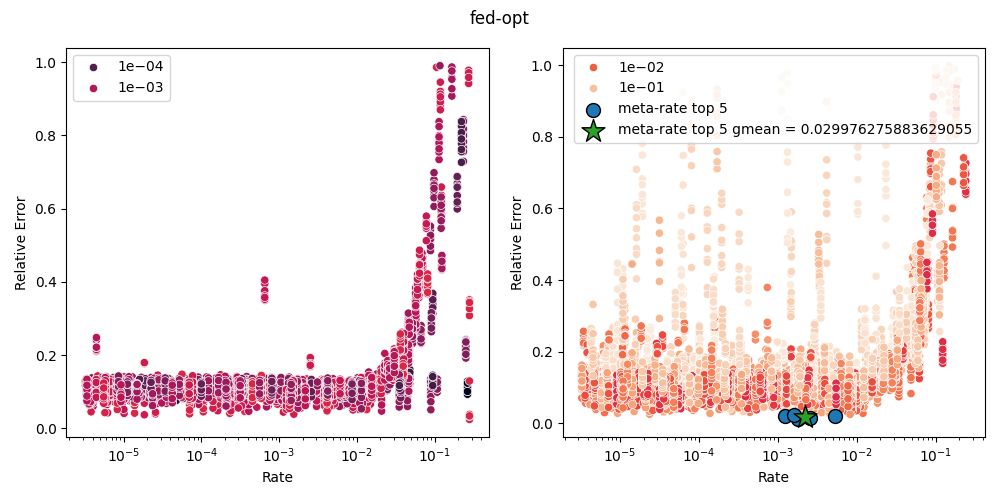}
    \caption{$x$: \textit{rate}; \textit{hue}: \textit{meta-rate}}
    \label{fig:experiments:dakota-poisson:fed-opt:rate-meta-rate}
  \end{subfigure}\\
  \begin{subfigure}[t]{0.85\textwidth}
    \centering
    \includegraphics[width=\textwidth]{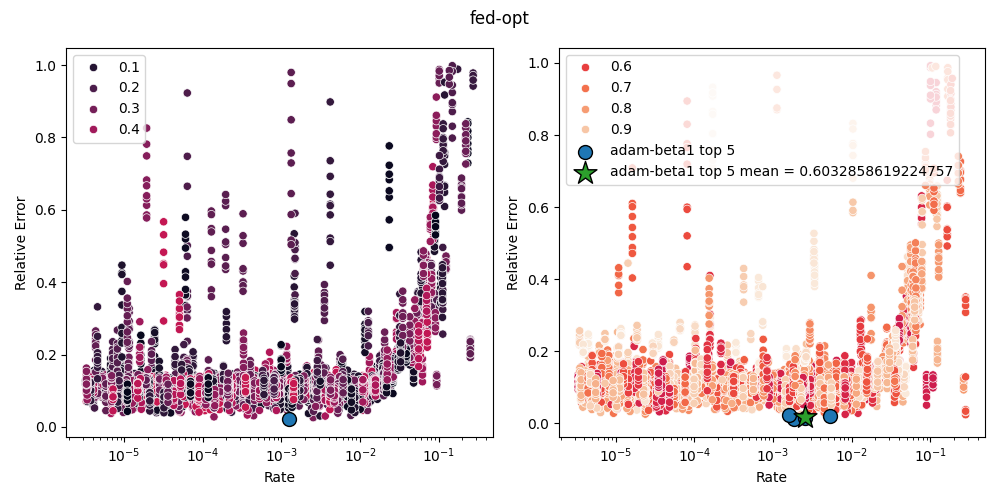}
    \caption{$x$: \textit{rate}; \textit{hue}: \textit{adam-beta1}}
    \label{fig:experiments:dakota-poisson:fed-opt:rate-adam-beta1}
  \end{subfigure}\\
  \begin{subfigure}[t]{0.85\textwidth}
    \centering
    \includegraphics[width=\textwidth]{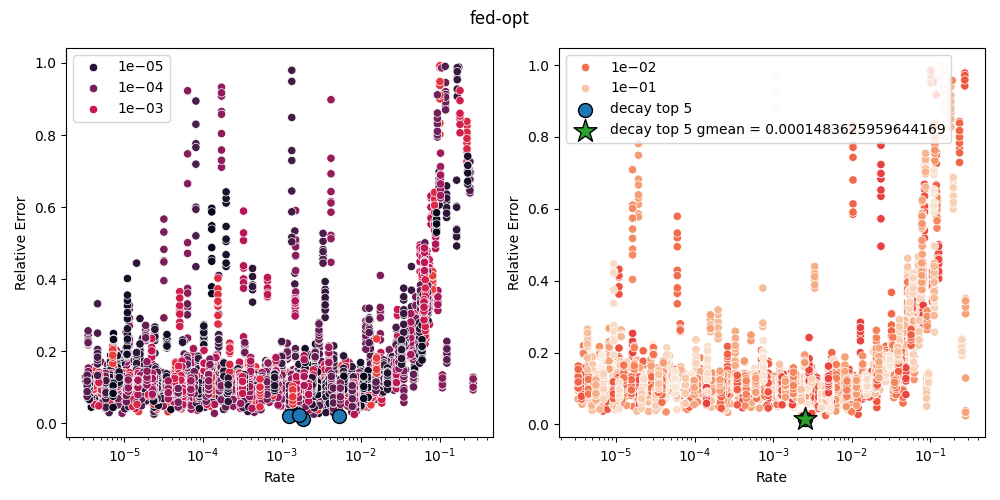}
    \caption{$x$: \textit{rate}; \textit{hue}: \textit{decay}}
    \label{fig:experiments:dakota-poisson:fed-opt:rate-decay}
  \end{subfigure}
  \caption{Scatter plots varying two parameters for {\gcpfed}. The $y$-axis
    denotes the relative error from the MLE; the first parameter of interest is
    plotted along the $x$-axis; variations in the second parameter of
    interest are
  denoted by the hue.}
  \label{fig:experiments:dakota-poisson:fed-opt}
\end{figure}

It is important to note that small correlation coefficient does not necessarily
mean a lack of dependence. The scatter plots
in~\cref{fig:experiments:dakota-poisson:gcp-sgd:rate-adam-beta1}--\cref{fig:experiments:dakota-poisson:fed-opt}
help identify trends that the correlation coefficients might miss. These plots
illustrate several multi-way relationships. The LHS sample values for the first
parameter, e.g., \textit{rate}, run along the $x$-axis. The $y$-axis corresponds
to an output response, i.e., relative error from the MLE or wallclock time. The
hue of each point in the scatter plot corresponds to a second parameter, e.g.,
\textit{decay}. We emphasize the data points corresponding to the top-5 settings
and their mean/geometric mean in the second parameter.

For any given output response type, the data is split on the second parameter so
that the data in the left plot corresponds to samples where the second parameter
values are below the median and conversely on the right. In several cases,
splitting the data in this way helped identify interesting patterns that were
missed by the correlation coefficients. For instance, the bifurcation
in~\cref{fig:experiments:dakota-poisson:gcp-sgd:rate-adam-beta1} suggests that
{\gcpsgd} is largely insensitive to variations in the \textit{rate}, which is
what one might interpret from the correlation coefficients
in~\cref{fig:experiments:correlations:gcp-sgd}. On the other hand, the
bifurcation indicates that the algorithm is more sensitive to
\textit{adam-beta1} than the correlation coefficients suggests: below
a threshold, convergence to the MLE is uncertain. A similar bifurcation appears
in~\cref{fig:experiments:dakota-poisson:gcp-sgd:rate-decay}, showing sensitivity
in the \textit{decay}. The timing results in these two plots indicate strong
nonlinear correlation in the lower portion of the splits and weaker correlation
in the upper portion, the latter of which corresponds to the range of parameters
where convergence to the MLE is more probable.
\Cref{fig:experiments:dakota-poisson:fed-opt} shows similar results for
{\gcpfed}, except splitting on the median sample value for each parameter is not
particularly illuminating. Taking their means, we see that the sets of
parameters that optimize each algorithm are different. These settings are
summarized in~\cref{tab:lhs:poisson:tuned}.

\begin{table}[!ht]
  \centering
  \caption{Tuned parameter values.}
  \label{tab:lhs:poisson:tuned}
  \begin{tabular}{lrr}
    \toprule
    Parameter & {\gcpsgd} & {\gcpfed} \\
    \midrule
    \textit{rate}       & $0.035817$                  & $0.002185$ \\
    \textit{decay}      & $0.123013$                  & $0.000148$ \\
    \textit{adam-beta1} & $0.630219$                  & $0.603286$ \\
    \textit{adam-beta2} & $0.966996$                  & $0.991181$ \\
    \textit{adam-eps}   & $1.388768 \times 10^{-13}$  & $5.523877
    \times 10^{-14}$ \\
    \textit{meta-rate}  & -                           & $0.029976$  \\
    \bottomrule
  \end{tabular}\\
  \textsuperscript{\textdagger} {\gcpfed} only
\end{table}

\subsubsection{Experiments using tuned hyperparameters}
To assess the generality of the relationship between the hyperparameters, degree
of asynchrony, and convergence to the MLE, we computed decompositions for both
algorithms using the synthetic count dataset and two real-world datasets taken
from the FROSTT repository~\cite{frosttdataset},
arXiv and amazon-reviews. These experiments employed both
default and tuned hyperparameters and varying degrees of asynchrony and
distributed parallelism. The combination of these three factors are the
principal object of this study.

\begin{figure}
  \centering
  \begin{subfigure}[t]{0.9\textwidth}
    \centering
    \includegraphics[width=\textwidth]{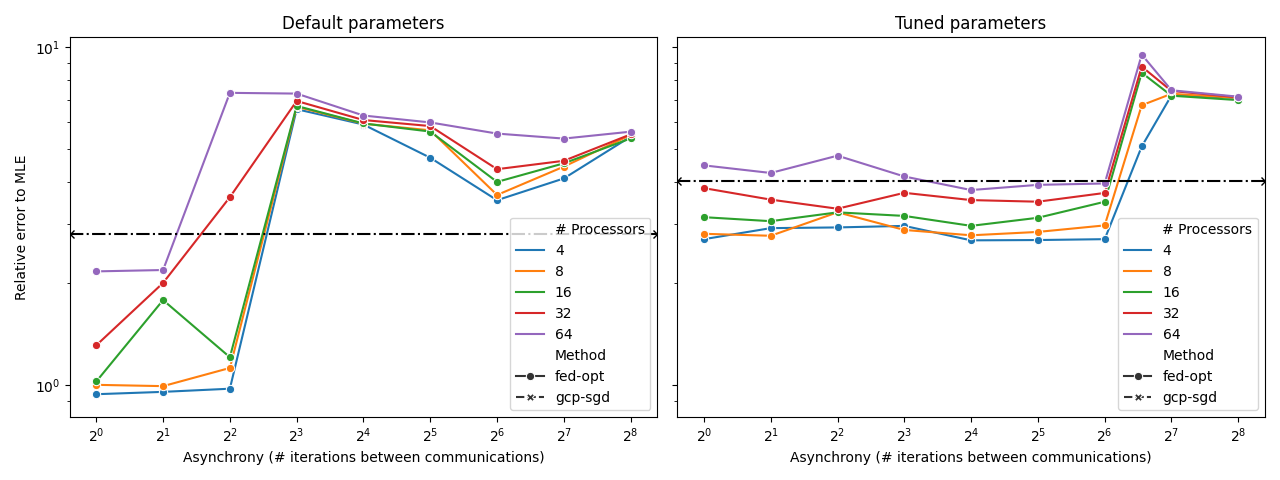}
    \caption{Relative error}
    \label{fig:experiments:dakota-poisson-valid:rel-err}
  \end{subfigure}\\
  \begin{subfigure}[t]{0.9\textwidth}
    \centering
    \includegraphics[width=\textwidth]{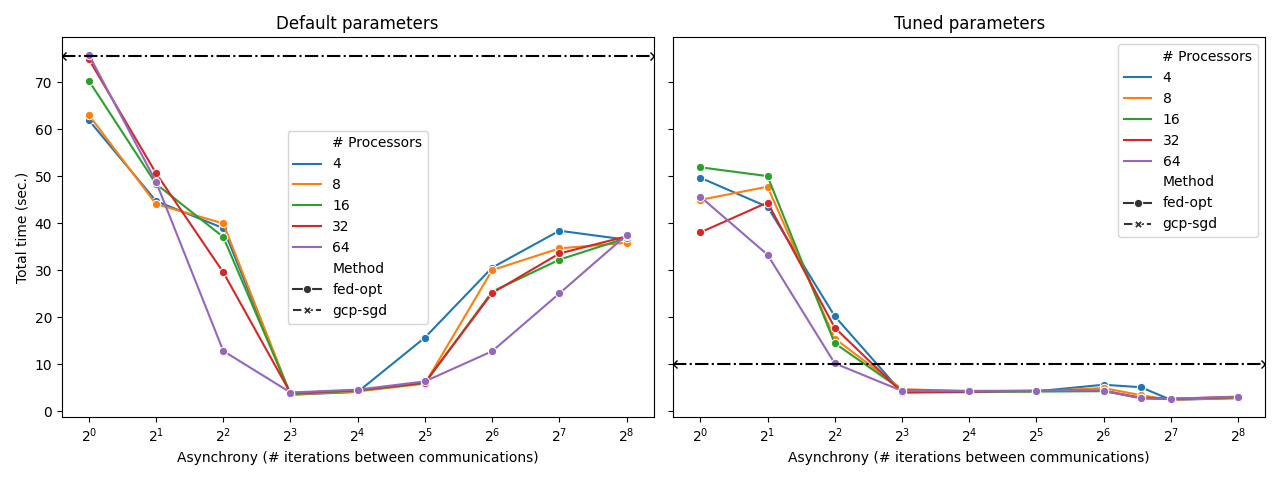}
    \caption{Total time}
    \label{fig:experiments:dakota-poisson-valid:time}
  \end{subfigure}
  \caption{Results on the synthetic count dataset
    comparing {\gcpsgd} and {\gcpfed} while varying the degree of asynchrony.
    Results include trials with both default and tuned parameters taken
  from the initial LHS study where only the random initial guess is different.}
  \label{fig:experiments:dakota-poisson-valid}
\end{figure}

\begin{figure}
  \centering
  \includegraphics[width=\textwidth]{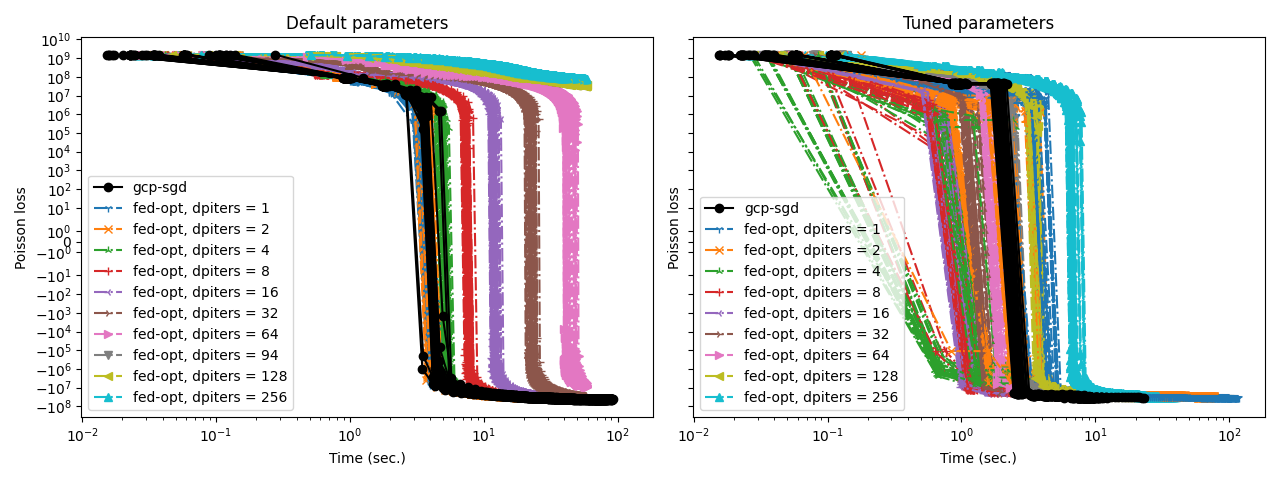}
  \caption{Traces of Poisson loss values by time in seconds for the
  arXiv dataset. The hue, \textit{dpiters}, denotes {\gcpfed} asynchrony.}
  \label{fig:experiments:arxiv:time-x-loss}
\end{figure}

\begin{figure}
  \begin{subfigure}[t]{0.9\textwidth}
    \centering
    \includegraphics[width=\textwidth]{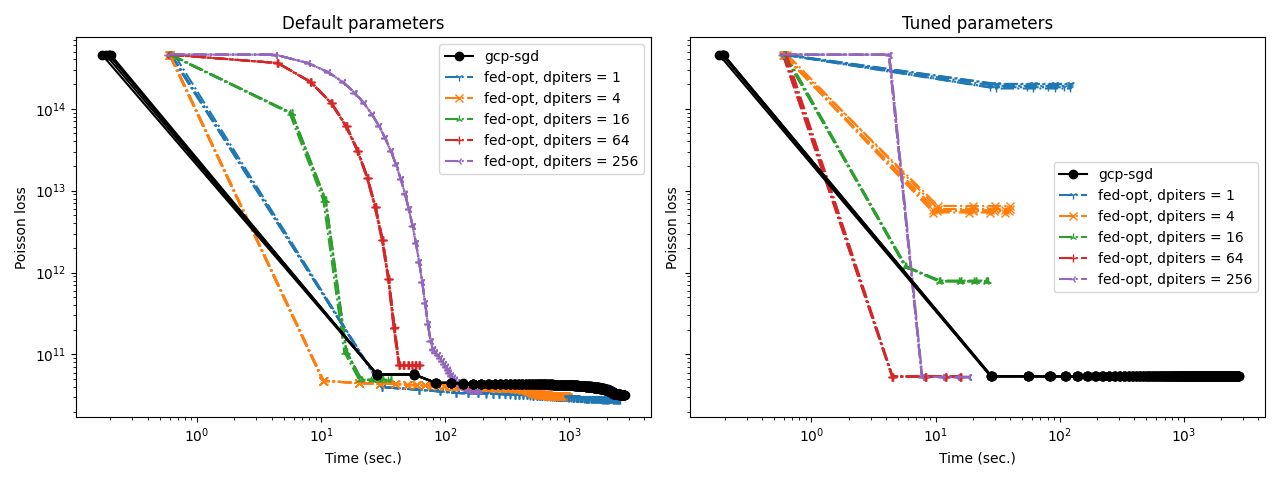}
    \caption{Traces of Poisson loss values by time in seconds.}
    \label{fig:experiments:amazon:time-x-loss}
  \end{subfigure}\\
  \begin{subfigure}[t]{0.9\textwidth}
    \centering
    \includegraphics[width=\textwidth]{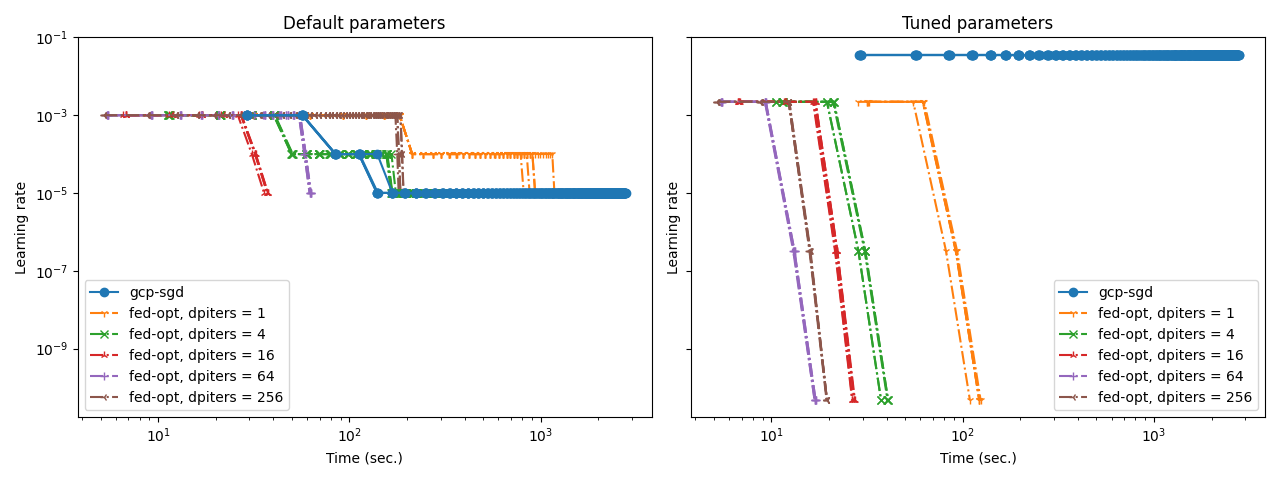}
    \caption{Traces of the {\adam} \textit{rate} parameter.}
    \label{fig:experiments:amazon:step}
  \end{subfigure}
  \caption{Results on the amazon-reviews dataset comparing {\gcpsgd} and
    {\gcpfed} while varying the degree of asynchrony.
    \Cref{fig:experiments:amazon:time-x-loss} compare the Poisson
    loss value against
    time. \Cref{fig:experiments:amazon:step} show the changes in {\adam}
    \textit{rate} in each iteration. Results include trials with both
  default and tuned parameter settings.}
  \label{fig:experiments:amazon}
\end{figure}

The trade-offs between convergence and cost as a function of hyperparameter
tuning, asynchrony, and distributed parallelism are illustrated for several
datasets in
\Cref{fig:experiments:dakota-poisson-valid}--\Cref{fig:experiments:amazon}.
Results with {\gcpfed} are presented as solid lines that vary along the $x$-axis
by asynchrony and in color by distributed parallelism. Results with
{\gcpsgd} are presented as horizontal dash-dot lines since it is synchronous.

In our analysis of the synthetic count dataset, we found both
algorithms converge
to two minima: the MLE and some other local minimum very far from the
MLE. \Cref{fig:experiments:dakota-poisson-valid} displays results
averaged from $n=10$ initial guesses. The low relative error
(\cref{fig:experiments:dakota-poisson-valid:rel-err}, left)
corresponds to trials where {\gcpfed} converged near to the MLE. We
found that the regime with higher relative error corresponds to
trials where {\gcpfed} terminated far from the MLE and hit the
maximum number of iterations. This suggests that infrequent
communication between clients and server required a greater number of
corrective epochs to reconcile local progress. This is expected since
the synthetic generation process, adapted from Chi and Kolda~\cite[\S
6.1]{Chi12TensorsSparsityNonnegative}, scales random factor matrix
entries to be much larger than nearby entries to generate the input
tensor; a synthetic count tensor would have greater data homogeneity
across factors otherwise.

When the constraint on accuracy can be relaxed, we see that {\gcpfed}
can be up to 21$\times$ faster than {\gcpsgd}. However, parameter
tuning seems to have a deleterious effect on accuracy for both
methods overall, which is surprising given all that changed was the
random initialization.

\paragraph{arXiv} 
Next, we consider the arXiv dataset. We performed a similar
comparison as above. GenTen was run with both default and tuned
parameters starting from $n=10$ initial guesses. We varied the
asynchrony and the number of processors. Each node in the cluster
consists of one 44-core IBM Power9 processor and 4 NVIDIA V100 GPUs
with 32 GB of RAM each. We compiled
GenTen with NVCC 11.2.152 and IBM Spectrum MPI
10.3.1.03rtm0 with Kokkos CUDA back-end. We ran GenTen MPI+CUDA with
4 MPI processes per node, each bound to a V100, with the number of
nodes ranging from $N = 1, 2, \ldots, 16$.

The results in~\Cref{fig:experiments:arxiv:time-x-loss} show the
traces in Poisson loss versus time isolating for asynchrony (i.e.,
{\gcpfed} \textit{downpour-iters}) for all $N$ nodes, since there was
no demonstrable effect for change in the number of processors. These
results demonstrate a regime where asynchrony and the
finite-difference global model update is advantageous
\textit{vis-a-vis} synchronous {\gcpsgd}. The comparable convergence
and performance behavior with and without tuned hyperparameters
suggests the data is fairly homogeneous across processors. In both
situations, default or tuned parameters, it is possible that the
asynchronous method can find solutions comparable to the synchronous
method. Unlike in the case with default parameters, it is interesting
that greater asynchrony results in worse computational performance
for the tuned parameters. However, there is a limit to the benefit of
asynchrony, as {\gcpfed} is unable to find the basin of attraction of
the MLE when sychronization is very infrequent.


\paragraph{Amazon}
We conclude our convergence study by considering the Amazon dataset.
We performed these experiments distributed across the same 18
compute nodes as above. Due to its large size, we utilized the
MPI+CUDA configuration only, to achieve reasonable solution times for
$N=5$ initial guesses with different parameterizations.

\Cref{fig:experiments:amazon} presents the traces of the Poisson loss
value and {\adam} \textit{rate} parameter for each
parameterization/initial guess combination. The behaviors on the left
(default parameters) are expected: both algorithms exit when they
fail to reduce the loss after 3 attempts. Thus {\gcpsgd} and
{\gcpfed} terminate with the same \textit{rate}. For low asynchrony
(i.e., $\textit{downpour-iters}=1$), {\gcpfed} convergence closely
resembles {\gcpsgd}. This is expected, due to the finite-differences
computed to reconcile local factors globally by {\gcpfed} in every
iteration. As in the case with the arXiv dataset, while asynchrony
increases, the impact of these pseudo-gradients emerges: {\gcpfed}
requires more computation to find the basin of attraction of the MLE
since it must reconcile progress made locally that is adverserial to
the global model. This suggests this dataset exhibits increased data
heterogeneity across processors.

The behaviors with parameters tuned (right) for another problem
provide a cautionary tale. Here, asynchrony has the benefit of
avoiding convergence to non-MLE local minima and results in up to a
172$\times$ speedup over synchronous {\gcpsgd}. It is clear
from~\cref{fig:experiments:amazon} (right) that greater asynchrony
interacts positively with lower initial learning rate. It may be
reasonable to conclude that asynchrony can be beneficial in the
convergence accuracy-performance trade-off in situations when the
initial learning rate is untuned.

%% file: conclusions.tex
In this work, we described on-node and distributed parallelism for {\gcpsgd}
implemented in GenTen. We provided performance results for several
configurations on real-world datasets. We also introduced {\gcpfed}, a federated
learning approach for {\gcpsgd}. We explored the role of asynchrony and
pseudo-gradient estimates in distributed environments and leveraged Latin
hypercube sampling to better understand its trade-offs. In our specific cases,
we found that hyperparameter tuning is typically problem-dependent and may not
generalize to problems with the same data distributional assumptions.
Additionally, despite tuning parameters for solution accuracy, we
found that improved computational performance was the principal beneficiary.
Despite this, we stress that the complex nature of computing CP tensor
decompositions---a nonlinear, nonconvex optimization problem---warrants novel
hyperparameter tuning that transcends techniques used in traditional machine
learning and uncertainty quantification for scientific simulation applications
alike.